\documentclass[12pt]{amsart}
\usepackage{amsmath,amsfonts,amssymb,amsthm,amsxtra,pstricks}
\usepackage[english]{babel}         

\usepackage{epsfig,epic,eepic,graphics,rotating,graphicx,bbm,texdraw,color}
\usepackage[all]{xy}               
\newcommand{\M}{{\mathbb M}}
\newcommand{\N}{{\mathbb N}}
\newcommand{\A}{{\mathbb A}}

\newcommand{\C}{{\mathbb C}}

\renewcommand{\P}{{\mathbb P}}
\newcommand{\D}{{\mathbb D}}

\newcommand{\E}{{\mathbb E}}

\newcommand{\bm}{{\mathbf M}}
\newcommand{\bt}{{\mathbf T}}
\newcommand{\bx}{{\mathbf X}}

\newcommand{\be}{{\mathbf E}}
\newcommand{\bbf}{{\mathbf F}}
\newcommand{\bu}{{\mathbf U}}
\newcommand{\bw}{{\mathbf W}}

\newcommand{\kf}{{\mathcal F}}

\newcommand{\ki}{{\mathcal I}}

\newcommand{\kl}{{\mathcal L}}
\newcommand{\km}{{\mathcal M}}
\newcommand{\kn}{{\mathcal N}}
\newcommand{\ko}{{\mathcal O}}

\newcommand{\kq}{{\mathcal Q}}

\newcommand{\ku}{{\mathcal U}}

\newcommand{\kz}{{\mathcal Z}}

\newcommand\coker{\operatorname{coker}\nolimits}

\DeclareMathOperator{\id}{id}

\DeclareMathOperator{\Supp}{Supp}

\DeclareMathOperator{\Hom}{Hom}
\DeclareMathOperator{\Hilb}{Hilb}

\DeclareMathOperator{\Spec}{Spec}

\DeclareMathOperator{\length}{length}
\DeclareMathOperator{\Grass}{Grass}
\DeclareMathOperator{\Sing}{Sing}

\newcommand{\ttop}{T_{top}}

\newcommand{\lra}{\longrightarrow}

\newcommand{\xra}{\xrightarrow}
\newcommand{\isoto}{{\lra\hspace{-1.5 em}
\raisebox{ 0.6 ex}{$\textstyle\sim$}\hspace{0.8 em}}}
\newcommand{\intoo}[1]{\:
\xymatrix@1{\ar@{^(->}[r]^{#1}&}\:
}

\begin{document}
\newtheorem{sub}{}[section]
\newtheorem{subsub}{}[sub]
\newtheorem{subsubsub}{}[sub]

\title[Bubble tree compactification]{bubble tree compactification of moduli
spaces of vector bundles on surfaces}

\author[Markushevich]{D.~Markushevich}
\address{Math\'ematiques - b\^{a}t.M2, Universit\'e Lille 1,
\newline F-59655 Villeneuve d'Ascq Cedex, France}
\email{markushe@math.univ-lille1.fr}

\author[Tikhomirov]{A.S.~Tikhomirov}
\address{Department of Mathematics, State Pedagogical University,
\newline Respublikanskaya Str 108
\newline 150 000 Yaroslavl, Russia}
\email{astikhomirov.mail.ru}

\author[Trautmann]{G.~Trautmann}
\address{Universit\"at Kaiserslautern, Fachbereich Mathematik,
\newline  Erwin-Schr\"odinger-Stra{\ss}e
\newline D-67663 Kaiserslautern}
\email{trm@mathematik.uni-kl.de}

\footnotetext[1]{The first author acknowledges the support of the French Agence Nationale de
Recherche VHSMOD-2009 Nr. ANR-09-BLAN-0104. Research of the second and third author was supported by the
DFG Schwerpunktprogramm 1094}

\begin{abstract}
In this article we announce some results on compactifying moduli spaces
of rank-2 vector bundles on surfaces by spaces of vector bundles on trees of
surfaces. This is thought as an algebraic counterpart of the so called
bubbling of vector bundles and connections in differential geometry.
The new moduli spaces are algebraic spaces arising as quotients
by group actions according to a result of Koll\'ar.
As an example the compactification of the space of stable rank-2 vector
bundles with Chern classes $c_1=0, c_1=2$ on the projective plane is studied
in more detail. Proofs are only indicated and will appear in separate
papers.
\end{abstract}

\vspace{2cm}
\maketitle

\tableofcontents

\begin{section}{Introduction}

In this article, we describe a conceptual scheme of a new
construction of a compactification of moduli spaces of stable bundles on
surfaces whose boundary consists of vector bundles on trees of surfaces,
replacing the torsion free semistable sheaves appearing in the Gieseker-Maruyama compactification.
In this description the long proofs of completeness, separatedness,
versality, properness are replaced by brief sketches, and the complete versions
will appear in full detail in subsequent papers. In conclusion, we produce a concrete
example of the compactification of the moduli space of stable bundles on
the projective plane $\P_2$ with second Chern class $c_2=2$.
In this example, we provide an alternative explicit construction of the same compactification and
give a complete description of its boundary.

To some extent, the replacement of limit sheaves in a compactification
by vector bundles on trees of bubbles is very natural.
The bubbling phenomenon appeared in eighties and nineties in the description
of degeneration processes in several conformally invariant
problems of geometric analysis: minimal surfaces (Sacks--Uhlenbeck),
harmonic maps (Parker),
pseudoholomorphic curves in symplectic varieties (Parker--Wolfson, Rugang Ye),
and Yang--Mills fields on 4-manifolds (Feehan \cite{Fe}, Taubes \cite{Ta1}, Uhlenbeck \cite{U}).

Donaldson--Uhlenbeck constructed a (partial) compactification $\overline{YM_n}$ of the moduli
space $YM_n$ of instantons of charge $n$ on a 4-manifold $S$, where the instantons are defined as
the $ASD$ Yang--Mills connections on a vector
bundle over a 4-manifold $S$ with Chern classes $c_1=0$,
$c_2=n$:
\[
\overline{YM_n}\subset YM_n\cup (YM_{n-1}\times X) \cup (YM_{n-2}\times S^2X)\cup
\cdots\cup S^nX,
\]
In \cite{MT1} and \cite{MT2} it was shown that
 in the case of $S=S^4$,  $\overline{YM_n}$ has a real semi-algebraic
structure.

The boundary of $\overline{YM_n}$ consists of ``ideal instantons'', that is, singular connections
whose curvature is a sum of a smooth part and of several delta-functions.
A degenerating family of bundles and connections on them encloses more information than is
kept by an ``ideal instanton''.
The bubble-tree compactification of Feehan--Taubes--Uhlenbeck (FTU) is a kind of blowup of the boundary of $\overline{YM_n}$, encoding the way an ASD connection degenerates into an ideal instanton
by a connection over a tree of surfaces, obtained from $S$ by successive gluings of spheres $S^4$ at a finite set of points.

When $S$ is a complex projective surface, Donaldson
proved that the Kobayashi--Hitchin correspondence identifies $YM_n$ with the 
moduli space of $\mu$-stable vector bundles $M^{\mu\textrm{-}s}(0,n)$ with 
$c_1=0$,  $c_2=n$ on $S$ (\cite{D1}; see also \cite{LT}
for further developments).  Thus the natural question can be asked, whether
the Uhlenbeck--Donaldson and the FTU compactifications also have an 
algebro-geometric interpretation.
For the first one, the answer is known:
Jun Li \cite{JL} endowed $\overline{YM_n}$ with a structure of a 
quasi-projective scheme
and defined a birational morphism
$\overline{M^{\mu\textrm{-}s}(0,n)}\to\overline{YM_n}$,
where the closure of $M^{\mu\textrm{-}s}(0,n)$ is taken in the 
Gieseker--Maruyama moduli space of semistable sheaves on $S$. See also 
\cite{HL}, Sect. 8.2, where a similar compactification, called moduli space of
$\mu$-semistable sheaves, is constructed for the sheaves with arbitrary 
determinant.

The main motivation of our work is to find an algebro-geometric analog of the 
FTU
compactification, which would be a kind of blowup of 
$\overline{M^{\mu\textrm{-}s}(0,n)}$.
By analogy with the topological bubbles which are 4-spheres, we can introduce
the notion of an algebraic bubble. Given a complex surface $S$ and its blowup
$\widetilde{S}$
in a point $p$ with the exceptional line $L$, an algebraic bubble is a complex 
projective plane
$P=\P_2(\C)$, attached to $\widetilde{S}$ along $L$.
Algebraic bubble trees are obtained by iterating this construction. They 
appear as fibers of the
semi-universal family over the compactified configuration spaces of 
Fulton--MacPherson \cite{FM}.

In the realm of algebraic geometry,
the degenerations are described in terms of
flat families.
It turns out that using flat deformations, one can ``replace'' singular 
sheaves $\kf$ on $S$ by
bundles on the trees of surfaces $S_T$, and that
the thus obtained tree bundles together with bundles on the original surface 
$S$
fit into a separated algebraic space of finite type. This is the main result 
of the paper, stated as Theorem
\ref{thmexist}. In some particular cases, for example, if $S=\P_2$, we can 
assert
that this algebraic space is proper (Theorem \ref{thmexist-a}).

This result still does not contain an answer to the question which served
its motivation, by the following reasons. First, our construction does not 
provide any morphism between our moduli space
and the Gieseker--Maruyama compactification in any direction. Second, though 
one can associate a topological bubble tree with 4-spheres as bubbles to any 
algebraic bubble tree by contracting to points the intersection lines of the 
algebraic bubbles, this does not lead to a correspondence
between the bundles on algebraic bubble trees
and the ASD connections on the associated spherical bubble trees.
The problem is that the stock of tree-like bundles we admit in our 
compactification
contains bundles which are nontrivial on
the intersection lines of bubbles, so there is no way to push them forward
along the topological contraction in order to get a bundle on a spherical
bubble tree. This is certainly a drawback of our approach based on the Serre 
construction
for flat families of rank 2 bundles over a curve.
It is a challenging problem to find another approach which would
bring us to a stock of boundary bundles, trivial on the intersections of 
components.

We will briefly mention some related work.
Tree-like bundles were used with the opposite goal
by D.Gieseker in \cite{G} in order to construct bundles on $S$
by deforming bundles on trees to bundles on $S$. Over curves, Nagaraj and 
Seshadri \cite{NS}
considered bundles on a degenerating family of curves and compactified them
by bundles on reducible curves, pasting in trees of rational components.
Buchdahl \cite{B1}, \cite{B2} studied, by differential geometric means, a 
compactification
of degenerating bundles on $S$ by bundles on a blowup of points in $S$ 
(which is an irreducible surface, unlike our bubble trees).

In this article, $S$ is always a smooth complex projective
surface, endowed with a very ample polarization class $h$, and
$M_{S,h}(2;\mathcal{N},n)$ denotes the Gieseker-Maruyama
moduli space of semistable torsion free sheaves $E$ of
rank 2 on $S$ with fixed Chern classes $c_1( E)=\mathcal{N}$ and
$c_2(E)=n\ge(\mathcal{N}^2)/4$, where ``semistable''
means ``Gieseker semistable with respect to $h$''. The space we are 
compactifying is the open subspace
$M^s_{S,h}(2;\mathcal{N},n)$ of stable vector bundles,
assumed nonempty.

If not otherwise stated, all the {\em schemes} are locally of finite type 
over $\C$, and the base
of any family we consider is always assumed to be a {\em scheme}.
\end{section}

\begin{section}{Trees of surfaces}

\begin{sub}\label{trs}{\bf Trees.}\rm\ 
A {\bf tree} $T$ in this article is a finite graph, oriented by a partial
order $<$ and satisfying:
\begin{itemize}
\item there is a unique minimal vertex $\alpha\in T$, the root of $T$;
\item for any $a\in T,\; a\neq \alpha$, there is a unique maximal vertex
$b<a$, the predecessor of $a$, denoted by $a^-$;
\item By $a^+:=\{b\in T\ |\ b^- =a\}$ we denote the set of direct
successors of $a\in T$. We let $\ttop$ denote the vertices of $T$
without successor.
\end{itemize}

A {\it weighted tree} is a pair $(T, c)$  of a tree $T$ with a map $c$
which assigns to each vertex $a\in T$ an integer $n_a$, called
the {\it weight} or {\it charge} of the vertex, subject to the conditions
\begin{equation}\label{trfirst}
n_a\geq 0 \text{ if } a\ne\alpha,
\end{equation}
\begin{equation}\label{trscond}
\#a^+\geq 2\ \text{ if }\ n_a=0 \text{ and } a\ne\alpha,
\end{equation}
\begin{equation}\label{trthird}
n_\alpha\geq C,
\end{equation}
where $C\leq0$ is some constant depending on $S$, $h$ and $\kn$, as specified below in formula
(\ref{constC}).
The total weight or total charge of a weighted tree is the sum
\;$\Sigma_{a\in T}\; n_a=n$ of all the weights.
We denote by $\bt_n$ the set of all trees which admit
a weighting of total charge $n$. It is obviously finite.

\end{sub}

\begin{sub}\label{trsurf}{\bf Trees of surfaces.}\rm\ 
Let $S$ be a smooth complex projective surface with an ample invertible sheaf
$\ko_S(h)$. Our trees of surfaces are reduced
surfaces $S_T$, defined for any tree $T$ and whose components are
indexed by the vertices of the tree $T$. They are constructed by the
following data.

\begin{itemize}
\item[(i)] For each vertex $a\not\in\ttop\cup \{\alpha \}$, let $P_a$ be a copy of $\P_2(\C)$
together with a line $l_a\subset P_a$ and a finite subset
$Z_a=\{x_b\in P_a\smallsetminus l_a|\; b\in a^+\}$, and let
$S_a\xra{\sigma_a} P_a$ be the blowup of $P_a$ along $Z_a$. We will denote the
exceptional lines in $S_a$ by $\tilde{l_b},\; b\in a^+,$
and the inverse image of $l_a$ in $S_a$ by the same symbol $l_a$.
\item[(ii)] For $a=\alpha$, we set $S_\alpha\xra{\sigma_\alpha} S$ to be
the blowup of $S$ at a finite set $Z_\alpha =\{x_a\in S|\; a\in \alpha^+\}.$
\item[(iii)] For $a\in\ttop$, $S_a$ is a copy of $\P_2(\C)$.
\item[(iv)] For each  $a> \alpha $, we fix some isomorphism
$\tilde{l_a}\xra{\phi_a} l_a$.
\end{itemize}

A tree-like surface of type $T$  over $S$ or a
$T$-surface is now defined as the result $S_T$ of gluing the above surfaces
$S_a$ along the isomorphisms $\phi_a$
into a reduced connected normal crossing surface
with components $S_a$. We write
$$S_T=\cup_{a\in T} S_a$$
and identify $\tilde{l_a}$ with $l_a$ for all $a\in T$.
After this identification the lines $l_a$ are the intersections
$l_a=S_a\cap S_{a^-}.$

By the construction of $S_T$, all or a part of its components
can be contracted. In particular, there is the morphism
\begin{equation}\label{sigmaS}  S_T\xra{\sigma} S\end{equation}
which contracts all the components except $S_\alpha$ to the points of the finite set $Z_\alpha.$

Note that:\\
1) There are no intersections of the components other than the lines $l_a$.\\
2) If $a\in\ttop$ then $S_a$ is a plane $\P_2(\C)$.\\
3) If $T=\{\alpha\}$ is trivial, then $S_T=S.$\\
4) After contracting the lines $l_a$ topologically (over $\C$), one obtains a
tree of $4$-sphere bubbles.
\end{sub}

\begin{sub}\label{linebu}{\bf Line bundles on $T$-surfaces.}\rm\ 
Let now $L$ be a line bundle on a $T$-surface $S_T$, and let the
inverse image of the divisor class $h$ on $S_\alpha$ also be denoted by $h$.
Then the restrictions of $L$ to the
components can be written as
 \[
L|S_\alpha=\ko_{S_\alpha}(mh-\Sigma_{a\in
  \alpha^+}\; m_a l_a)
\]
and
\[
L|S_a=\ko_{S_a}(m_al_a-\Sigma_{b\in a^+}\; m_b l_b).
\]

$L$ is called ample if each $L|S_a$ is ample for all $a\in T$. This means that
\[
m, m_a>0 \text{ and } m^2>\Sigma_{a\in\alpha^+}\; m^2_a
\text{ and } m^2_a>\Sigma_{b\in a^+}\; m^2_b.
\]
for all $a\in T.$
Let $m,r$ be positive integers. We will say that $L$ is
of type $(r,m,h)$ if its restriction to the root surface is given by
$L|S_\alpha=\ko_{S_\alpha}(rmh-r\Sigma_{a\in\alpha^+}\; m_a l_a) $ for some $m_a$
$(a\neq \alpha)$.

We define the multitype of a line bundle on a T-surface $S_T$ to be the sequence
$$ m_T:=(m_a)_{ a\in T},$$ where $m_{\alpha}=m$
The above inequalities imply
\end{sub}

\begin{sub}\label{lem1}{\bf Lemma.}\;
\begin{enumerate}
\item For any $m,r>0$ and any weighted tree $(T,c), T\in\bt_n$,
there are at most finitely many ample line bundles of type $(r,m,h)$.
\item Given $n$, there is an integer $m_0$ such that for any $m\geq m_0$
and any $(T,c), T\in\bt_n,$  there is an ample line bundle on $S_T$
of type $(1,m,h)$.
\end{enumerate}
\end{sub}
\end{section}

\begin{section}{Tree-like bundles}\label{trlb}

\begin{sub}\label{defvb}{\bf Definition.}\rm\quad
Let $S_T=\cup_{a\in T} S_a$ be a $T$-surface with tree
$T\in\bt_n$. A vector bundle $E=E_T$ on $S_T$ or the pair $(E_T,S_T)$ is
called a
tree of vector bundles or a $T$-bundle if the restrictions
$E_a=E|S_a,\, a\in T$, satisfy the following conditions:
\begin{itemize}
\item[(i)] $E_a$ has rank 2, $c_1(E_a)=0$ for $a\ne\alpha$,
$c_1( E_\alpha)=\sigma^*_\alpha\mathcal{N}$, and the second Chern classes
$c_2(E_a)= n_a,\; a\in T,$ define a weighting of $T$ in the sense of
definition \ref{trs}, (\ref{trscond}).
\item[(ii)] If $T\ne\{\alpha\}$, then $E_a$ is {\em admissible} as defined
below for any $a\in T.$
\item[(iii)] If $T=\{\alpha\}$, then $[E]\in M_{S,h}(2;\mathcal{N},n)$.
\end{itemize}
\end{sub}

\begin{sub}\label{hulsb}{\bf Definition.}\rm\quad
Let $S_T=\cup_{a\in T} S_a$ be a $T$-surface with tree
$T\in\bt_n$. Assume that $T\ne\{\alpha\}.$ Let  $a\in T$,
and let $E_a$ be a rank-2 vector bundle on $S_a$.
We will say that  $E_a$ is
{\it admissible} if one of the following conditions is satisfied:

\begin{itemize}
\item[(i)]
In case $a\in\ttop$, with $S_a$ a plane $P_a$,
$E_a$ is an extension of type
\begin{equation}\label{ext1}
0\to\ko_{P_a}\to E_a\to\ki_{x,P_a}\to0,\; x=\{pt\}\not\in l_a,\;
c_2(E_a)=1,
\end{equation}
\item[(ii)] or in case $a\in\ttop$, with $S_a$ a plane $P_a$,
$E_a$ is an extension of type
\begin{equation}\label{ext2}
0\to\ko_{P_a}(-1)\to E_a\to\ki_{Z,P_a}(1)\to 0,
\end{equation}
where $\dim Z=0,\; Z\cap l_a=\emptyset,\; c_2(E_a)=\length(Z)-1\ge 2.$
\item[(iii)]
In case $a\not\in\ttop$, $a\ne\alpha$, $|a^+|\ge 2$, $E_a$ is a non-split extension
of type
\begin{equation}\label{ext3}
0\to\ko_{S_a}(-l_a + \underset{b\in a^+}\Sigma l_b)\to E_{S_a}\to
\ki_{Z, S_a}(l_a-\underset{b\in a^+}\Sigma l_b)\to 0,
\end{equation}
where $\dim Z\le0,\; Z\subset S_a\smallsetminus\{(\underset{b\in a^+}\cup l_b)
\cup l_a\},\; c_2(E_{S_a})=\length(Z)+|a^+|-1\ge 1,$
\item[(iv)]
or $a\not\in\ttop$, $a\ne\alpha$, and $E_a$ is a non-split extension of the type
\begin{equation}\label{ext4}
0\to\ko_{S_a}(\underset{b\in a^+}\Sigma l_b)\to E_a\to
\ko_{S_a}(-\underset{b\in a^+}\Sigma l_b)\to0,\ \ \
c_2(E_a)=|a^+|\ge 1,
\end{equation}
\item[(v)]
or $a\not\in\ttop$, $a\ne\alpha$, and $E_a=2\ko_{S_a}$.\vskip3mm
\item[(vi)]
In case $a=\alpha$, $E_\alpha$ is a non-split extension of type
\begin{equation}\label{ext5}
0\to\ko_{S_\alpha}(-qh + \underset{b\in\alpha_0^+}\Sigma l_b)\to E_\alpha\to
\ki_{Z, S_\alpha}(qh-\underset{b\in\alpha_0^+}\Sigma l_b+\sigma_\alpha^*\mathcal{N})\to 0,
\end{equation}
\begin{equation}\label{constC}
c_2(E_\alpha)\ge C:=-q_0^2(h^2)-q_0|(h\cdot\mathcal{N})|,
\end{equation}
for some subset $\alpha_0^+$ of $\alpha^+$,
where $0\leq q\leq q_0$ for some integer $q_0$ depending on $S$,
and where $\dim Z\le0,\; Z\subset S_\alpha\smallsetminus(\underset{b\in\alpha^+}\cup l_b)$,
$\length(Z)\le n+q_0^2(h^2)+q_0|(h\cdot\mathcal{N})|$.
\end{itemize}
\vskip3mm

{\bf Remark.}\; The above definitions single out a possibly redundant class of vector bundles
including all the bundles which may
occur in degenerations. The conditions (i)--(vi) guarantee the boundedness of
the family of tree-like bundles and replace the (semi)stability conditions,
which are not obvious for tree-like bundles. Those $T$-bundles which indeed occur in degenerations will be called limit $T$-bundles, see Definition \ref{limbun}.
\vskip3mm

{\bf Notation.} Let $E=E_T$ be a $T$-bundle on $S_T$ for a tree $T\in\bt_n$.
We also write $E_T=\underset{a\in T}\# E_a$
and define its total second Chern class to be
$$c_2(E_T):= \underset{a\in T}\Sigma c_2(E_a)= \underset{a\in T}\Sigma n_a=n.$$

Two $T$-bundles $(E_T, S_T)$ and $(E'_T, S'_T)$ are called
isomorphic if there exists an isomorphism $\phi : S_T\to S'_T$ over $S$
such that $E_T = \phi^* E'_T$.
\vskip3mm

In the following we need formulas for the Euler characteristics
of line bundles and $T$-bundles. These follow by standard computations
on the components of the trees.
\end{sub}
\vskip5mm

\begin{sub}\label{chi}{\bf Lemma.}\;
For any $T$-bundle $(E_T, S_T)$, $T\in\bt_n,$ and any
ample line bundle $L$ on $S_T$ of type $(1,m,h)$, the following Euler
characteristics are independent of the tree and are given by the formulas
\begin{equation}\label{chi1}
\chi(E\otimes L)=m^2(h^2)+m(h\cdot(2\mathcal{N}-K_S))+\frac{1}{2}(\mathcal{N}\cdot(\mathcal{N}-K_S))+
\end{equation}
$$2\chi(\ko_S)-n=:N_m$$
\begin{equation}\label{chi3}
\chi(L)=\frac{m^2}{2}(h^2)-\frac m2(h\cdot K_S)+\chi(\ko_S)
\end{equation}
\end{sub}
\vskip5mm

Later we will consider the embeddings of $S_T$ into the Grassmannian $G=\Grass(N_m,2)$ of 2-dimensional
quotient spaces of the space $\C^{N_m}$, defined by the global sections of  an appropriate twist
of tree-like bundles $E_T$ on $S_T$. Here $N_m$ denotes the integer given by \eqref{chi1}. The universal rank-2 quotient bundle
on $G$ will be denoted by $\kq$.
We have the following boundedness result.
\vskip5mm

\begin{sub}\label{boundn}{\bf Proposition}\; {\rm (Boundedness)}.
For any $n$ there is an integer $m_0>0$ such
that, for any $m\geq m_0$ and any $T\in\bt_n$, there is an ample line bundle
$L$ of type $(1,m,h)$ on $S_T$ such that for any $T$-bundle $E_T$ on $S_T$,
\begin{enumerate}
\item [(i)] $h^i(E_T\otimes L)=0$ for $i>0$ and
$h^0(E_T\otimes L)=\chi(E_T\otimes L)=N_m$;
\item [(ii)] the evaluation map
$\C^{N_m}\otimes \ko_{S_T}\twoheadrightarrow E_T\otimes L$ is surjective;
\item [(iii)] the induced map $S_T\xrightarrow{i_G} G=\text{Gr}(\C^{N_m}, 2)$
is a closed embedding such that
$$i_G^\ast\kq\simeq E_T\otimes L\ ,\ i_G^\ast\ko_G(1)\simeq L^{\otimes 2}
\otimes\mathcal{O}_{S_T}(\sigma^*\mathcal{N});$$
\item [(iv)] for any $j,q>0$,\, $h^j(i_G^*\ko_G(q))=0,$ and the Hilbert
polynomial $P_G(q):=\chi(i_G^*\ko_G(q))$ is given by the formula
\begin{equation}\label{P_G}
P_G(q)=2q^2m^2(h^2)+2q^2m(h\cdot\mathcal{N})-qm(h\cdot K_S)
\end{equation}
$$+\frac{1}{2}q^2(\mathcal{N}^2)-\frac{1}{2}q(\mathcal{N}\cdot K_S)+\chi(\mathcal{O}_S(\mathcal{N})).$$
\end{enumerate}

\begin{proof} This follows from Serre's theorems A and B,
the boundedness of the family of all
admissible bundles of given type on each $S_a$, which follows easily from the definition,
and from the boundedness of the family of $h$-semistable vector bundles \cite{Sim}
with $c_1=\mathcal{N}$ and $c_2\le n$ on $S.$
(ii) follows from the above and \cite[lemma 5.13]{T}.
Formula (\ref{P_G}) follows directly
from \ref{chi}.
\end{proof}

\end{sub}
\end{section}
\vskip1cm

\section{Families of tree-like bundles}\label{treefam}

In this section we fix the definition of families of $T$-surfaces $S_T$
and $T$-bundles  $(E_T,S_T)$ for trees $T\in\bt_n$ with fixed
total charge $n$.


\begin{sub}\label{goodfam}{\bf Definition}\rm\;
($\bt_n$-families of surfaces).
1) Let $\mathbf{X}\overset{\pi}{\to} Y$ be a flat family of trees of surfaces
over a scheme $Y$ locally of finite type,
whose trees belong to ${\bf T}_n$. Such a family is called a
family of trees of surfaces of type ${\bf T}_n$, or simply a
{\it $\bt_n$-family},
if there exists a morphism $\sigma:\mathbf{X}\to S\times Y$
such that the following holds:

(i) $\pi=pr_2\circ\sigma$.

(ii) For each closed point $y\in Y$ the morphism
$\sigma_y=\sigma|S_y:S_y\to S\times\{y\}\simeq S$,
where $S_y=\pi^{-1}(y)$, is the standard contraction (\ref{sigmaS}),

(iii) There is a union of irreducible components $\bx^b$ of $\bx$
such that the restriction $\sigma|\bx^b$ is a birational morphism
$\bx^b\to S\times Y.$

The morphism $\sigma:\mathbf{X}\to S\times Y$ will also be called
{\it a standard contraction}.



2) A $\bt_n$-family of surfaces $\mathbf{X}\overset{\pi}{\to} Y$ is called a
{\it good} $\bt_n$-family of surfaces if $\sigma$ is birational on the whole
of $\bx$.


3) A $\bt_n$-family of surfaces $\mathbf{X}\overset{\pi}{\to} Y$ is called
trivial if $\sigma$ is an isomorphism.

4) Two $\bt_n$-families
$\mathbf{X}\overset{\pi}{\to}Y$ and
$\mathbf{X'}\overset{\pi'}{\to}Y$ over the same base $Y$
are called {\it isomorphic} if there exists an isomorphism
$\mathbf{X}\xra{\phi}\mathbf{X'}$ such that
\begin{equation}\label{rigidity}
\pi'\circ\phi=\pi\quad and \quad\sigma'\circ\phi=\sigma,
\end{equation}
where $\mathbf{X}\overset{\sigma}{\to}S\times Y$ and
$\mathbf{X'}\overset{\sigma'}{\to}S\times Y$ are
the standard contractions.
\end{sub}
\vskip3mm

\begin{sub}\label{explfam}{\bf Examples.}\rm\; 1) The standard contraction
$\mathbf X=S_T\to S$ of a single $T$-surface for a tree $T\in\bt_n$ is a $\bt_n$-family over a point with
$\mathbf X^b= S_\alpha$, $\sigma|\bx^b =\sigma_\alpha:S_\alpha\to S$.
This $\bt_n$-family is good only if
$T=\{\alpha\}$.

2) Let $C$ be a smooth curve and let $X\to S\times C$ be the blowup
of a point $(s,c)$. Then  $X\to C$ is a good $\bt_1$-family.

3) Let $\bx=S_T\times Y$ be the product of a $T$-surface with some scheme $Y$ 
locally of finite type,
$T\in\bt_n$, $T\ne\{\alpha\}$.
Then $\bx\to Y$ is a $\bt_n$-family  with
$\mathbf X^b= S_\alpha\times Y$ which is not good.
\end{sub}
\vskip2mm
For families over smooth curves we have the following

\begin{sub}\label{aksing}{\bf Theorem.}\;
Let $\bx\to C$ be a good $\bt_n$-family of surfaces over a smooth
curve $C$. Then $\bx$ has at most $A_k$-singularities, analytically locally 
trivial along the lines of intersection
of components in the fibres of $\bx$.

{\rm The $A_k$-singularities in $\bt_n$-families over curves really appear in 
the
construction of limit bundles in the proof of the Completeness Theorem
\ref{complthm} as a result of certain contractions of tree-like surfaces.}
\end{sub}
\vskip3mm

\begin{sub}\label{bdlfam}{\bf Definition}\rm\ \;
($\bt_n$-families of bundles).
(a) For a fixed total second Chern class $n$, a $\bt_n$-family of tree-like 
bundles is
given as a triple $(\be/\bx/Y)$, where $\bx/Y$ denotes a $\bt_n$-family
$\bx\overset{\pi}{\to} Y$ of surfaces and $\be$
is a vector bundle on $\bx$, such that
its restriction to all the components of all the fibers of $\pi$ over the 
closed points of $Y$ are admissible.

(b) Let now $\bx\overset{\pi}{\to} Y$ be a good $\bt_n$-family of surfaces. 
By Definition \ref{goodfam}, (ii), there exists a maximal dense
open subset $U$ of $Y$ with the property that the fibers $X_y=\pi^{-1}(y)$ are 
isomorphic to $S$ for all closed points $y\in U$.
A  triple $(\be/\bx/Y)$ as in (a) will be called a
{\it good $\bt_n$-family of (tree-like) bundles} if the restrictions of $\be$ 
to the fibers $X_y$ over the closed
points $y\in U$ are {\em stable} vector bundles from $M^s_{S,h}(2;\kn,n)$.

There is an obvious notion of isomorphy and equivalence for $\bt_n$-families.
Two $\bt_n$-families of bundles $(\be/\bx/Y)$ and $(\be'/\bx'/Y)$ are called 
isomorphic,
\[
(\be/\bx/Y)\simeq (\be'/\bx'/Y),
\]
if there is an isomorphism $\bx\xra{\phi}\bx'$ of $\bt_n$-families of
surfaces such that
$\be\simeq\phi^\ast\be'.$ The families are called equivalent if there is an
isomorphism $\bx\xra{\phi}\bx'$ and an invertible sheaf $L$ on $Y$ such that
$(\be\otimes\pi^*L/\bx/Y)$ and $(\be'/\bx'/Y)$ are isomorphic.
\end{sub}

\begin{sub}\label{limbun}{\bf Definition}\;\rm\  (Limit bundles).
A $T$-bundle $(E_T, S_T)$ of total charge $n$ is called a {\em limit} 
$T$-bundle if there exists
a (germ of a) smooth pointed curve $(C,0)$ and a good $\bt_n$-family
of bundles $(\be/\bx/C)$ as defined in \ref{bdlfam} (b), such that
$\bx|C\smallsetminus \{0\}\xra{\sigma}S\times (C\smallsetminus \{0\})$ is an
isomorphism, $\be|C\smallsetminus \{0\}$ is a family of stable vector bundles
from $M^s_{S,h}(2;\kn,n)$ and $E_T\simeq\be|\bx_0$ on the fibre of $\bx$
over $0\in C.$

A limit $T$-bundle is called {\em sss-limit} $T$-bundle if there exists
a (germ of a) smooth pointed curve $(C,0)$ and a $\bt_n$-family
of bundles $(\be/\bx/C)$ as defined in \ref{bdlfam} (a), such that
$\bx|C\smallsetminus \{0\}\xra{\sigma}S\times (C\smallsetminus \{0\})$ is an
isomorphism, $\be|C\smallsetminus \{0\}$ is a family of strictly semi-stable 
vector bundles
from $M_{S,h}(2;\kn,n)$ and $E_T\simeq\be|\bx_0$ on the fibre of $\bx$
over $0\in C.$

Let $\M_n(pt)$ denote the set of all limit $T$-bundles with
$T\in\bt_n$ and $\bm_n(pt)$ the set of their isomorphism classes. More generally,
we give the following definition.

\end{sub}

\begin{sub}\label{stack}{\bf Definition}\;\rm\  (Moduli stack and moduli functor).
For any $Y$ as above, we denote by $\M_n(Y)$ the set of all $\bt_n$-families
$(\be/\bx/Y)$ of limit bundles. Given a morphism $f: Y^\prime\to Y$,
it is obvious how to define the pullback $f^\ast(\be/\bx/Y)$ of a family
$(\be/\bx/Y)$ and it is easy to verify that this is again a $\bt_n$-family
of limit bundles.
Thus $\M_n$ is a pseudofunctor
$\M_n: (Sch/\C)^{op}\to (Sets)$  in the language of stacks, where in our
setting
$(Sch/\C)$ denotes the category of complex schemes locally of finite
type over $\C$.
By definition, $\M_n(pt)$ can be identified with $\M_n(\Spec(\C))$.\\
To get a functor, we define $\bm_n(Y):= \M_n(Y)/\sim$, where $\sim$
denotes equivalence.
\end{sub}

\begin{sub}\label{subfunctors}{\bf Open sub-pseudo-functors of $\M_n$.}\rm\
Consider the open sub-pseudo-functors $\M_n^g$, $\M_n^t$, $\M_n^0$ and $\M_n^s$ of the pseudo-functor
$\M_n: (Sch/\C)^{op}\to (Sets)$, defined as follows:
$$
\M_n^g(Y)=\{(\be/\bx/Y)\in\M_n(Y)\ |\ \forall\ y\in Y,\ \be_y\ \mbox{is not an sss-limit
bundle}\},
$$
$$
\M_n^t(Y)=\{(\be/\bx/Y)\in\M_n(Y)\ |\ \bx\to Y\ \text{is a trivial}\ \bt_n\text{-family}\},
$$
$$
\M_n^0(Y)=\M_n^g(Y)\cup \M_n^t(Y),
$$
$$
\M_n^s(Y)=\{(\be/\bx/Y)\in\M_n(Y)\ |\ \bx\to Y\ \text{is a trivial}\ \bt_n\text{-family, and}
$$
$$
\ \forall\ y\in Y,\ \be_y\in M^s_{S,h}(2;\kn,n)\}\ \ \ \ (Y\in Sch/\C).
$$
Obviously, for any $Y\in Sch/\C$,
\begin{equation}\label{intersn of psf}
\M_n^s(Y)=\M_n^g(Y)\cap\M_n^t(Y).
\end{equation}

The pseudofunctors $\M_n^g$, $\M_n^t$, $\M_n^0$ and $\M_n^s$ and their associated functors $\bm_n^g$, \ldots will be used
in the construction of the moduli spaces $M_n^g$, $M_n^0$, see Section \ref{modsp}.

\end{sub}

\begin{sub}\label{lift}{\bf Lift of families.}\;\rm\ 
In the proofs of the main results one mostly has to deal with good families of
tree-like bundles over curves.
We use the standard notation $(C,0)$ for a curve $C$ with a marked
point $0\in C$ and denote by $C^*$ the punctured curve $C\smallsetminus\{0\}$.
By a finite covering of curves $\tau:(\tilde C,0)\to (C,0)$
we understand a finite morphism
$\tau:\tilde C\to C$ such that $\tau(0)=0$.

For a given curve $(C,0)$ and a given family of tree-like bundles
$\bbf=(\be/\bx/C)\in\M_n(C)$ over $C$
we denote the {\it restriction of $\mathbf{F}$ onto $C^*$} by
$$\bbf^*=(\be^*/\bx^*/C^*)\in\M_n(C^*),$$
where
$\bx^*=C^*\times_C\bx,\, \pi^*=\pi|\bx^*,\, \sigma^*=\sigma|\bx^*,\,
\be^*=\be|\bx^*.$

Respectively, we denote the lift of $\bbf=(\be/\bx/C)\in\M_n(C)$ to
$\tilde{C}$ by
$$\tilde{\bbf}=\tilde{C}\times_C\bbf=(\tilde{\be}/\tilde{\bx}/\tilde{C})
\in\M_n(\tilde{C}),$$ where
$\tilde{\bx}=\tilde{C}\times_C\bx,$ and
$\bx\overset{\lambda}\leftarrow\tilde{\bx}
\overset{\tilde{\pi}}\rightarrow\tilde{C}$
denote the natural projections with
$\tilde{\be}=\lambda^*\be$.

Moreover, the notation for the {\it lift} $\tilde{\bbf}\xra{\tilde{\phi}}
\tilde{\bbf'}$ to $\tilde{C}$ of an isomorphism $\bbf\xra{\phi}\bbf'$
over $C$ should be selfexplaining.
\end{sub}

The proofs of properness and separatedness of the moduli space
we are going to construct are reduced to the respective properties
of the pseudofunctor $\M_n$ over the smooth curves,
basing upon the valuative criteria for properness and separatedness.
Thus the following completeness and separatedness theorems for the pseudofunctor
$\M_n$ are key results for the existence of a
proper moduli space of tree-like bundles.
The proofs can only roughly be indicated in this note.

\begin{sub}\label{complthm}{\bf Completeness Theorem.}\;
The pseudo-functor $\M_n$ is complete in the following sense.
For any smooth pointed curve $(C,0)$ and any family of tree-like bundles
$(\be/C^*\times S/C^*)\in \M_n(C^*)$,
there is a finite covering
\[
(\widetilde{C}, 0)\xrightarrow{\tau} (C,0)
\]
and a family $(\widetilde{\be}/\widetilde{\bx}/\widetilde{C})\in
\M_n(\widetilde{C})$ of tree-like bundles
together with an isomorphism $\varphi$,
\[
\xymatrix{
C\times S\ar[d] & \widetilde{\bx}\ar[l]\ar[d]\ar @{-^{)}} @<-2.8pt> [r]
\ar @{-} @<2.2pt> [r]
 &
\widetilde{\bx}^\ast\ar[d]\ar[r]^-\varphi_-\approx & \widetilde{C}^*\times
S\ar[dl]\\
C & \widetilde{C}\ar[l] \ar @{-^{)}} @<-2.8pt> [r]
\ar @{-} @<2.2pt> [r]&
\widetilde{C}^\ast &
}
\]
such that
$$\widetilde{\be}|\widetilde{X}^\ast\cong \varphi^*(\tau\times
\id_S)^*\be.$$
\end{sub}
\vskip5mm

{\em Sketch of the proof.}
\textbf{Step 1}: Let $\bx=C\times S$. We use the description 
of reflexive rank-2 sheaves $\bbf$ on $\bx$ via the
Serre construction\footnote{For $S=\P_2$ one might use monads.} in the relative
situation:
\[
0\to\ko_\bx(-qh)\xrightarrow{s}\bbf\to \ki_Z(qh+\mathcal{N})\to 0,
\]
where $Z$ is finite over $C$.

Consider the closure $\overline{{M}_{S,h}^s(2;\mathcal{N},n)}$ of the moduli
space of stable vector bundles
${M}_{S,h}^s(2;\mathcal{N},n)$ in the
projective scheme of Gieseker--Maruyama of $\chi$-semistable torsion-free 
sheaves. Using
its construction and projectivity, we can replace $C$ by a finite covering and
assume that there is a reflexive sheaf $\be$ over the whole of $X$ such that:
\begin{itemize}
\item $\be|C^\ast\times S$ is the given bundle;
\item there is a number $q$ with an exact sequence
\[
0\to \ko_\bx(-qh)\to \be\to\ki_Z(qh+\mathcal{N})\to 0\:;
\]
\item $Z$ is reduced, smooth over $C^\ast$, and  $Z\to C$ is flat and finite.
\end{itemize}
\vskip5mm

\textbf{Step 2}: Smoothing $Z$ by blowups of points of
$Z\cap(\{0\}\times S)$ in the $3$-fold $\bx=C\times S$.
\vskip5mm

\textbf{Step 3}: Separating the components of $Z$ by lifting the families
over finite base changes and by further blowing up,
thus getting ``bubbles'' and separating the corresponding extension
sequences of the Serre construction.
\vskip5mm

\textbf{Step 4}: Contracting superfluous bubbles in the $3$-fold obtained in Step 3, thereby producing a 1-parameter family of tree-like surfaces whose total space is a three-dimensional variety having at worst curves of $A_k$ singularities.

\textbf{Step 5}: The case of a $\bt_n$-family $(\be/\bx^*/C^*)$ such that all the fibres of $(\bx^*/C^*)$ are reducible. Since all the bundles
$E_y,\ y\in C,$ are limit bundles, after a possible shrinking of $C^*$, there exists a surface $Y_0$ containing
$C^*$ and a $\bt_n$-bundle $(\be_0/\bx_{Y_0}/Y_0)\in\M_n(Y_0)$ such that, for $Y^*=Y_0\smallsetminus C^*$ one has
$\bx_{Y_0}\times_{Y_0}Y^*\simeq S\times Y^*$. Extend $Y_0$ to a surface $Y$ containing $C$ and, respectively, extend
$(\be_0|S\times Y^*/S\times Y^*/Y^*)$ to a family of Gieseker semistable sheaves
$(\be_Y/S\times Y/Y)$. As in Step 1, represent $\be_Y$ as a sheaf obtained by the Serre construction, applied to some subscheme $Z_Y$ of $S\times Y$, finite over $Y$. Then take
$Z=Z_Y\times_YC$ and apply Steps 2--4 to the family $(\be_Y|S\times C/S\times C/C)$. As a result, possibly after a base change over $C$, we obtain a $\bt_n$-family $(\widetilde{\be}/\widetilde{\bx}/C)$ extending
$(\be/\bx^*/C^*)$.

\vskip5mm
\begin{sub}\label{septhm}{\bf Separatedness Theorem.}\;
The pseudofunctor $\M_n^g$ is separated in the following sense:
Let $(C,0)$ be a smooth pointed curve and let
$\bbf=(\be/\bx/C),\;\bbf'=(\be'/\bx'/C)\in\M_n^g(C)$
be two
families of tree-like bundles.
Suppose that there is an isomorphism
$$\xymatrix{\bx^*=C^*\times_C\bx\ar[dr]_{\ \ \ \sigma_1} \ar[rr]^{\phi_0}_{\sim}&&
 \ar[dl]^{\sigma_2\ \ \ }C^*\times_C\bx'=\bx{'^*}\\
&C^*\times S&}$$
such that the restrictions
$\bbf^*=(\be^*/\bx^*/C^*),\; \bbf^{'*}=(\be^{'*}/\bx^{'*}/C^*)\in\M_n^g(C^*)$
of $\bbf,\bbf'$ over $C^*$ are isomorphic via
$\phi_0$, i.e. there exists an isomorphism of vector bundles $\psi_0:\be^*\isoto\phi_0^*\be^{'*}$.
Then $\phi_0$ extends to an isomorphism $\phi:\bx\isoto \bx'$, and
there exists an invertible sheaf $L$ on ${C}$ with an isomorphism
$\rho:L|{C}^*\isoto\mathcal{O}_{{C}^*}$
such that the isomorphism $\psi_0\otimes\pi_0^*(\rho):\be^*\otimes \pi_0^*(L|{C}^*)\isoto\phi_0^*\be^{'*}$
extends to an isomorphism
$\be\otimes\pi^*L\isoto\phi^*\be'$, which provides an isomorphism
of the $\bt_n$-families of bundles
$(\be\otimes\pi^*L/\bx/C)\simeq (\be'/\bx'/C)$.
\end{sub}

The proof is rather elaborate, and what follows gives a brief idea of it.
The first step is blowing up $\bx$, $\bx'$ to obtain a model $\tilde\bx$,
smooth over $\C$ and dominating both $\bx$, $\bx'$. There is an isomorphism of the lifted
vector bundles $\mu_0:\tilde\be^*\isoto\tilde\be^{'*}$ over $C^*$.
We can twist $\be$ by $\pi^*(\ki^k)$ for some integer $k$, where $\ki$ is the ideal sheaf of $0\in C$,
so that $\mu_0$ extends to a sheaf morphism $\mu:\tilde\be\to\tilde\be'$. We assume that $k$ is chosen
to be minimal with this property, so that the restriction of $\mu$ to the fiber over $0\in C$ is
a nonzero morphism of sheaves $\mu_0:\tilde E_0\to \tilde E'_0$.

Our definitions imply that the determinants of  $\be$, $\be'$ are lifts of line bundles from $C\times S$,
and $\det\mu$ can be viewed as a
section of $(\det\tilde\be)^{-1}\otimes\det\tilde\be'\simeq\tilde\pi^*\kl$,
where $\kl$ is a line bundle on $C$.
This implies that the support of $\coker\mu$ is a simple normal crossing surface,
a union of components $S_a$ of $\tilde X_0$. A combinatorial argument shows
that $\det\mu$ cannot vanish only on a part of components,
so if $\mu$ is not an isomorphism, then $\mu_0$ is degenerate on every component of $\tilde X_0$.
It is quite obvious then that the image of $\mu_0$, restricted to $S_a$,
is a rank 1 sheaf over every component $S_a$ of $\tilde X_0$.

The second step is the proof of a fact from commutative algebra, which, stated in geometric terms, reads as follows:

\begin{sub}\label{commalg}{\bf Lemma.}\;
Let $\phi:E\to E'$ be an injective morphism of rank $2$ locally free sheaves on a smooth irreducible
threefold $X$,
and assume that $D=\Supp\coker\phi$ is an effective divisor on $X$ having smooth irreducible components
$D_i,\ i\in A$, such that $L_i:=\coker\phi|D_i$ is a rank $1$ sheaf for each $i\in A$. Then $L_i$
is a line bundle on $D_i,\ i\in A$. In particular, for each $i\in A$ the bundle $E'|D_i$ has a
quotient line bundle $L_i$.
\end{sub}

As follows from Definition \ref{hulsb} (i-ii), there is always a component in $\tilde X_0$, namely, the inverse image of any top component of $X'$, on which the restriction
of $\tilde E'$ has no invertible quotient. This proves that $\mu$ has to be an isomorphism.

The last step of the proof is an argument, showing that there is only one way to contract some of
the components in $\tilde X$, on which $\tilde E$ is trivial, in order that the result of the contraction
might satisfy the condition~\eqref{trscond}.

\begin{sub}\label{S-equiv}{\bf Remark.}\;\rm The separatedness fails for the full functor $\M_n$ as there may be non-isomorphic $S$-equivalent vector bundles on $S$ which are limits of the same family of stable ones.
\end{sub}

\section{Fulton--MacPherson configuration spaces}\label{config}

The construction of the moduli space, up to technical details,
follows the standard pattern.
First we use Hilbert schemes and embeddings into Grassmannians for constructing
a space $H$, which parametrizes all the objects we want to include in our moduli space,
and then we quotient $H$ by a group action. First, to construct the parameter space $H$, we invoke
the results of  W.\,Fulton and R.\,MacPherson from their paper
\cite[sect.1-3]{FM},
where they introduce the configuration spaces $S[\mathbf{n}]$ for any
natural number $\mathbf{n}.$

\begin{sub}\label{FMfam}{\bf Notation.}\;\rm
Let $\mathbf{n}:=\max\{|T|-1\ |\ T\in\mathbf{T}_n\}$, and let
$Y_{FM}=S[\mathbf{n}]$ be the Fulton-MacPherson configuration space with
the semi-universal family
$S[\mathbf{n}]^+$ over $S[\mathbf{n}]$ of the so-called
``$\mathbf{n}$-pointed stable configurations'' over $Y_{FM}$. Let
\begin{equation}\label{FM1}
\pi_{FM}:\mathbf{X}_{FM}:=S[\mathbf{n}]^+\overset{\sigma_{FM}}\lra S\times Y_{FM}\overset{pr_2}\to Y_{FM}.
\end{equation}
denote the standard morphisms,
where $\sigma_{FM}$ is a birational morphism which decomposes into a sequence
of blowups
with explicitly described smooth centers. The family (\ref{FM1}) has a
number of nice
properties. In particular, both $\mathbf{X}_{FM}$ and $Y_{FM}$ are smooth
projective varieties and the following ``versality'' property holds:
\end{sub}
\begin{sub}\label{versal}{\bf Proposition}\;{(Versality)}.
Let $X\to S\times Y\to Y$ be a flat deformation of a standard contraction
$S_T\to S$ with base point $y\in Y.$ Then there is an open neighbourhood
$V(y)\subset Y$ of the point $y$ and a morphism $V(y)\xra{f}Y_{FM}$ such that
$X_V\to S\times V(y)\to V(y)$ is the pull back of
$\bx_{FM}\to S\times Y_{FM}\to Y_{FM}$ under $f$.
\end{sub}
In particular, for any
$T\in\mathbf{T}_n$ and any tree of surfaces $S_T$
there exists a point (not unique) $y\in Y_{FM}$ such that
\begin{equation}\label{FM2}
S_y:=\pi_{FM}^{-1}(y)\simeq S_T.
\end{equation}

Note that $Y_{FM}$ represents a functor of configuration families,
\cite[Theorem 4]{FM}, but there is no obvious way of deriving the above
versality from that. A proof can be done by a detailed analysis of the
Kodaira-Spencer map related to this deformation problem.
\vskip4mm

The birational morphism
$\sigma_{FM}:\mathbf{X}_{FM}\to S\times Y_{FM}$ is by its construction
decomposed into a sequence
$\sigma_{FM}=\sigma_1\circ...\circ\sigma_R$
of blowups with smooth centers, say,
$Z_i\subset(\sigma_1\circ...\circ\sigma_{i-1})^{-1}(S_0\times Y_{FM}),\ i=1,...,R.$
By this construction, the divisors
$D_i:=(\sigma_i\circ...\circ\sigma_R)^{-1}(Z_i)$
satisfy the property that, for any $y\in Y_{FM}$,
$D_i\cap S_y$ is a subtree of the tree $S_y=\pi_{FM}^{-1}(y)$.
For a sequence of positive integers $m_0,n_1,...,n_R$ we define the
invertible sheaves
$$\km_0:=\ko_S(m_0h)\boxtimes\ko_{Y_{FM}}\; \text{on}\; S\times Y_{FM}$$
and
$$\km:=\sigma_{FM}^*(\km_0)\otimes\ko_{\bx_{FM}}(-\Sigma n_iD_i)\; \text{on}\;
\bx_{FM}.$$
Using the above property of the divisors and Serre's theorems A and B,
one can derive the following Lemma.

\begin{sub}\label{lemFM}{\bf Lemma.}\quad For the given number $n$ there
is a sequence of positive integers $m_0, n_1, \ldots, n_R$ such that
\begin{itemize}
\item[(i)] $\km_0^r$ is $pr_2$-very ample for any $r\geq 1$ and
$pr_{2\ast}\km_0^r$ is locally free.
\item[(ii)] $\km^r$ is $\pi_{FM}$-very ample for any $r\geq 1$ and
$\pi_{FM\ast}\km^r$ is locally free.
\end{itemize}
\end{sub}

\begin{sub}{\bf Remark.}\label{yy'}\;\rm
(i) Note that, since $Y_{FM}$ is a projective variety, then for each
$T\in\bt_n$ and each $a\in T$ the number
$\epsilon(T,a):=\max\{\epsilon\in\N|$\ there exists $y\in Y_{FM}$ such that
$S_y$ has $T$ as its graph and $m_a(\km|S_y)=\epsilon\}$ is clearly finite,
where, as above,
$S_y=\pi_{FM}^{-1}(y)$ and we use the notation from \ref{linebu}.

(ii) From the definition of the line bundle $\km$ it follows that, if
$y,y'\in Y_{FM}$ are two points such that the fibres as trees of surfaces
$S_y$ and $S_{y'}$ have the same graph $T$, then the line bundles
$\km|S_y$ and $\km|S_{y'}$ have the same multitype:
$\mathbf{m}_T(\km|S_y)=\mathbf{m}_T(\km|S_{y'})$. In particular,
$$S_y\simeq S_{y'}\quad \text{implies}\quad \km|S_y\simeq \km|S_{y'}.$$

We thus are led to the following notation for an arbitrary tree of
surfaces $S_T$ with $T\in\bt_n$:
\begin{equation}\label{multityp}
\mathbf{m}_T(\km):=\mathbf{m}_T(\km|S_y)
\end{equation}
for any isomorphism
$S_T\overset{\sim}\to S_y,\; y\in Y_{FM}$. This notation is coherent, for the right hand side
of \eqref{multityp} does not depend on the
choice of $y$.
\end{sub}

Next, since the set of all pairs $(E_T,S_T)\in\M_n(\C)$ is bounded by
Proposition \ref{boundn}, we may
strengthen the result of Lemma \ref{lemFM} in the following way.

\begin{sub}\label{propFM}{\bf Proposition.}\quad
One can choose the numbers $m_0,n_1,...,n_R$ in Lemma \ref{lemFM}
in such a way that, for any
 $(E_T,S_T)\in\M_n(\C)$ and any isomorphism
$\phi_y:S_T\overset{\sim}\lra S_y,\; y\in Y_{FM}$, the following holds:

\begin{itemize}\itemsep3mm
\item[(i)] $h^j(S_y, \km|S_y)=0$ for any $j>0$ and the number
$r_0:=h^0(S_y, \km|S_y)$ is independent of $y$ with the above property.
\item[(ii)] Put $m:=r_0m_0$. Then the line bundle $L:=\phi_y^*(\km^{r_0})$
on $S_T$ has type $(1,m,h)$,\; $h^j(E_T\otimes L)=0,\; j>0,$
and
$N_m:=h^0(E_T\otimes L)=\chi(E_T\otimes L)$ is given by
(\ref{chi1}).
\item[(iii)] For any isomorphism
$\theta_y:\C^{N_m}\isoto  H^0(E_T\otimes L)$,
the induced map
\begin{equation}\label{evalmap}
\theta(y):\C^{N_m}\otimes\ko_{S_T}\overset{\theta_y\otimes id}\lra
H^0(E_T\otimes L)\otimes\ko_{S_T}\overset{ev}\lra E_T\otimes L
\end{equation}
is surjective, and the induced morphism $i_{\theta(y)}:S_T\to G:=Grass(N_m,2)$ to
the Grassmannian of 2-dimensional quotients of $\C^{N_m}$ is an embedding
such that
\begin{equation}\label{pbQ}
i_{\theta(y)}^*\kq=E_T\otimes L\quad\text{and}\quad i_{\theta(y)}^*\ko_G(1)\simeq L^2
\otimes\mathcal{O}_{S_T}(\sigma^*\mathcal{N}),
\end{equation}
where $\kq$ is the universal rank 2 quotient sheaf on $G$.
\item[(iv)] $h^j(i_{\theta(y)}^*\ko_G(q))=0$ for all $j,q>0$\; and
\begin{equation}\label{}
\chi(i_{\theta(y)}^*\ko_G(q))=2q^2m^2(h^2)+2q^2m(h\cdot\mathcal{N})-qm(h\cdot K_S)
\end{equation}
$$+\frac{1}{2}q^2(\mathcal{N}^2)-\frac{1}{2}q(\mathcal{N}\cdot K_S)+\chi(\mathcal{O}_S(\mathcal{N}))$$
(the same value as in (\ref{P_G})).
\end{itemize}
\end{sub}
\begin{proof}
(i) follows by a standard argument and Lemma \ref{lemFM},
(ii) and (iv) follow directly from Lemma \ref{lemFM} and Riemann-Roch.
(iii) follows from a lemma on embeddings into Grassmannians in
\cite[lemma 5.13]{T}, and boundedness, Proposition \ref{boundn}.
\end{proof}

Note that from (\ref{pbQ}) it follows immediately that
\begin{equation}\label{multiG}
\mathbf{m}_T(i_{\theta(y)}^*\mathcal{O}_{G}(1)\otimes\ko_{S_T}(-\sigma^*\mathcal{N}))
=2\mathbf{m}_T(\phi_y^*\mathcal{M}_{FM}^{r_0}).
\end{equation}

\begin{sub}\label{canbu}{\bf A functorial line bundle.}\;\rm
By the above the sheaf $\pi_{FM\ast}\km$ is locally free on $Y_{FM}$ of rank
$r_0$.  We consider the line bundle
$$L_{FM}^c:=\km^{r_0}\otimes \pi_{FM}^*(det \pi_{FM\ast}\km)^{-1}$$
on $X_{FM}.$ Using this line bundle, one can construct for any $\bt_n$-family
$(\be/\bx/Y)$ a line bundle $L_Y$ on $\bx$ such that these line bundles
are compatible with base change in the sense of J.Koll\'ar,
\cite[Definition 2.3]{Ko}, see \ref{canbuH} and Section \ref{modsp}.
This defines a descent of the line bundles $L_Y$ to a line bundle $L_{M_n}$ on the moduli space.
One of the possible
approaches to the proof of the projectivity of our moduli space would be to show
that $L_{M_n}$ is ample. To this end, one might verify the
weak positivity property \cite{V} for the bundles $L_Y$, but this seems to be difficult for the $\bt_n$-families
$(\be/\bx/Y)$ that are not good.
\end{sub}
\vskip5mm

\section{The Hilbert scheme construction}
The parameter
space for $\bt_n$-bundles will be an open part of the Hilbert scheme
$\Hilb^{P_H}(S\times G)$ consisting of $\bt_n$-surfaces, where as above
$G$ is the Grassmannian $Gr(N_m,2)$, and $P_H$ will be determined below.

\begin{sub}\label{embtree}{\bf Definition.}\;\rm
Let $m=r_0m_0$ be as in \ref{propFM}. An embedded $\bt_n$-surface is defined
to be a closed embedding $S_T\overset{i}{\hookrightarrow} S\times G$ such that
\begin{itemize}
\item[(i)] the composition $S_T\overset{i_G}{\hookrightarrow} G$ is
a closed embedding with $i^*_G\ko_G(1)\otimes\ko_{S_T}(-\sigma^*\mathcal{N})$ ample of type $(2,m,h)$.
\item[(ii)] the composition $S_T\to S$ is a standard contraction
as defined in (\ref{sigmaS}), which is an isomorphism if $T$ is a single
vertex.
\end{itemize}
\end{sub}
\vskip3mm

Let $\ko_{S\times G}(1):=\ko_S(mh-\mathcal{N})\boxtimes\ko_G(1)$ be the chosen
very ample polarization of $S\times G$ and, for an embedded tree
$S_T\overset{i}{\hookrightarrow} S\times G$, let
$P_H(q)=\chi(i^*(\ko_{S\times G}(q)))$ be the corresponding Hilbert polynomial.
From Definition \ref{embtree} it follows immediately that
$i^*(\ko_{S\times G}(1))$ is a very ample line bundle of type $(3,m,h)$ on 
$S_T$.
Hence,
$P_H(q)$ is given by the formula (\ref{chi3}) with $qm$ substituted for $m$:
\begin{equation}\label{polyn}
P_H(q)=\frac{9}{2}q^2m^2(h^2)-\frac{3}{2}qm(h\cdot K_{S})+\chi(\mathcal{O}_{S}).
\end{equation}
Consider now the Hilbert scheme $\Hilb^{P_H}(S\times G)$ and let
$$H^{'}\subset\Hilb^{P_H}(S\times G)$$
be the open subscheme of all embedded  $\bt_n$-surfaces in the sense
of Definition \ref{embtree}, and let $H^s\subset H^{'}$ be the open part
of those $S_T\overset{i}{\hookrightarrow} S\times G$ for which $S_T\simeq S$ 
and $(i^*(\ko_S(-mh)\boxtimes Q)/S_T/\Spec(\C))\in\M_n^s(\Spec(\C))$.

Finally, let $\overline{H}^s$ be the closure of $H^s$ in $H'$ and define
$$
H:=\{(S_T\overset{i}\hookrightarrow S\times G)
\in\overline{H}^s\ |\ i^*(\ko_S(-\mathcal{N})\boxtimes\ko_G(1))=
\phi^*\km^{2r_0}\
$$
$$
\text{for an}\ \text{isomorphism}\ \phi:S_T\overset{\sim}\to S_y\ \text{for some point}\ y\in Y_{FM}\ \text{and}
$$
$$
(i^*(\ko_S(-mh)\boxtimes Q)/S_T/\Spec(\C))\in\M_n(\Spec(\C))\}.
$$
Here the existence of such a point $y\in Y_{FM}$ follows from the
versality of $Y_{FM}$ and the property that
$i^*(\ko_S(-\mathcal{N})\boxtimes\ko_G(1))\simeq\phi^*\km^{2r_0}$
does not depend on the choice of the point $y$ in view of (\ref{yy'}), (ii).
Equivalently,
$$
H=\{(S_T\hookrightarrow S\times G)
\in\overline{H}^s\ |\
\mathbf{m}_T(\ko_S(-\mathcal{N})\boxtimes\ko_G(1)|S_T)= 2\mathbf{m}_T(\km^{r_0})\ \text{and}
$$
$$
(i^*(\ko_S(-mh)\boxtimes Q)/S_T/\Spec(\C))\in\M_n(\Spec(\C))\}.
$$

Since $H$ is a locally closed subscheme of the Hilbert scheme, there is
the semi-universal
family $\bx_H$ of embedded $\bt_n$-surfaces with diagram
\begin{equation}\label{univX}
\xymatrix
{
{\bx_H}\ar[dr]_{\pi_H}\ar@{^{(}->}[r]^{i_H\ \ \ \ \ } & S\times G\times H\ar[d]^{pr_3}\\
& H.
}
\end{equation}

\begin{sub}\label{canbuH}{\bf Lemma.}\; There is a line bundle $L_H$ on
$\bx_H$ such that for any fibre
 $\bx_{H,z}$, $z=(S_T\hookrightarrow S\times G)$,
and any isomorphism $\phi: S_T\xra{\sim} S_y$ ($y\in Y_{FM}$),
$$L_H|\bx_{H,z}\simeq L^c_{FM}\simeq \km^{r_0}|S_y$$ and that the Hilbert
polynomials of the fibres $\bx_{H,z}$ are formed with respect to the
line bundle $L_H\otimes i_H^*(\ko_S(-\mathcal{N})\boxtimes\ko_G(1)\boxtimes\ko_H).$
\end{sub}

The bundle $L_H$ plays also the role of a functorial polarization
in the sense of J.Koll\'ar. It can be obtained as follows. By the versality
of $Y_{FM}$ there is an open covering $(H_i)$ of $H$ with morphisms
$H_i\xra{f_i}Y_{FM}$ such that $X_H|H_i$ is the pullback of $\bx_{FM}$
under $f_i$. Then the pulled back bundles $f_i^* L^c_{FM}$ can be glued
to give the bundle $L_H$.
Using the bundle $L_H$, we define the bundle
\begin{equation}\label{bundleEH}
\be_H:=\mathcal{O}_{S}\boxtimes Q\boxtimes\mathcal{O}_H|\bx_H
\otimes(L_H)^{-1},
\end{equation}
Then $(\be_H,\bx_H,H)$ is a $\bt_n$-family and belongs to $\M_n(H).$

\begin{sub}\label{boundn2}{\bf Remark}\;\rm (Boundedness). It follows from the
results of the next section that any $\bt_n$-bundle occurs in the family
$(\be_H/\bx_H/H)$, proving that the functor $\bm_n$ is bounded.
\end{sub}
\vskip5mm

\section{The coarse moduli space}\label{modsp}

Given an arbitrary family $(\be/\bx/Y)\in \M_n(Y)$, one can construct
a line bundle $L_Y$ on $\bx$ as in the case of the family over $H$ using
the versality of $Y_{FM}$. This bundle is fibrewise isomorphic to
$\km^{r_0}$ and has type $(1,m,h)$. If there is a morphism $\rho:Y^\prime\to Y$,
then the constructions of the line bundles $L_Y$ are compatible, so that
$L_{Y^\prime}\simeq \rho^*L_Y$.
By Proposition \ref{propFM},
$\pi_*(\be\otimes L_Y)$ is locally free of rank $N_m$, and we can consider
the principal $GL(N_m)$-bundle
$$\tilde{Y}:={\rm Isom}(\mathbf{k}^{N_m}\otimes\mathcal{O}_Y,
\pi_*(\mathbf{E}\otimes L_Y))\overset{\rho}\to Y$$
over $Y$ with the Cartesian diagram

\begin{equation}\label{cartes}
\xymatrix
{\tilde{\bx}\ar[d]_{\tilde{\pi}}\ar[r]_{\tilde{\rho}} & \bx\ar[d]_{\pi}\\
\tilde{Y}\ar[r]_{\rho} & Y}
\end{equation}

Denote the lifts of the bundles by $\mathbf{\tilde{E}}={\tilde{\rho}}^*\mathbf{E},\ L_{\tilde{Y}}=\tilde{\rho}^*L_Y$.
In view of Proposition \ref{propFM}, we have on $\mathbf{\tilde{X}}$ a
universal epimorphism
\begin{equation}\label{theta}
\Theta:\ \C^{N_m}\otimes\mathcal{O}_{\mathbf{\tilde{X}}}\overset{\sim}\to
\tilde{\pi}^*\tilde{\pi}_*(\mathbf{\tilde{E}}\otimes L_{\tilde{Y}})
\overset{ev}\twoheadrightarrow\mathbf{\tilde{E}}\otimes L_{\tilde{Y}},
\end{equation}
which induces an embedding
$i_\Theta:\mathbf{\tilde{X}}\hookrightarrow S\times G\times\tilde{Y}$
in the commutative diagram
\begin{equation}\label{hilpY}
\xymatrix
{
\mathbf{\tilde{X}}\ar[dr]_{\tilde{\pi}}\ar@{^{(}->}[r]^{i_\Theta\ \ \ \ \ } &
S\times G\times \tilde{Y}\ar[d]^{pr_3}\\
& \tilde{Y}
}
\end{equation}
such that
\begin{equation}\label{tildeE}
\mathbf{\tilde{E}}\otimes L_{\tilde{Y}}=i_\Theta^*(\ko_S\boxtimes Q\boxtimes
\ko_{\tilde{Y}})
\end{equation}
Let now $L_{\tilde{Y}}\otimes i_\Theta^*(\ko_S\boxtimes \ko_G(1)\boxtimes
\ko_{\tilde{Y}})$ serve as the polarization for computing the Hilbert 
polynomial
in $S\times G\times\tilde{Y}.$ The restriction of $L_{\tilde{Y}}$ to
each fibre of $\tilde{\bx}$ over a point $\tilde{y}\in\tilde{Y}$ is by its
construction isomorphic to some $\km^{r_0}_y\otimes\ko_{S_T}\otimes\ko_G(1)$
or to some $ \ko_{S_T}(mh)\otimes\ko_G(1)$. Hence the fibres of $\tilde{\bx}$
all have the Euler characteristic $P_H(q),$ and the
conditions of the definition of $H$ are satisfied for the fibers 
of $\tilde\bx$.
By the universal property of the Hilbert scheme, there is a morphism
$\phi$ in the following Cartesian diagram such that $\tilde{\bx}$ is the 
pull back
of $\bx_H$,
\begin{equation}\label{univmor}
\xymatrix
{\tilde{\bx}\ar[d]^{\tilde{\pi}}\ar[r]^{\tilde{\phi}} & {\bx_H}\ar[d]^{\pi_H}\\
\tilde{Y}\ar[r]^{\phi} & H}
\end{equation}
and such that
\begin{eqnarray}
\tilde{\be}\otimes L_{\tilde{Y}} & \simeq & i_\Theta^*(\ko_S\boxtimes Q
\boxtimes\ko_{\tilde{Y}})\nonumber\\  & \simeq & i_\Theta^*(id\times\phi)^*
(\ko_S\boxtimes Q\boxtimes\ko_H)\nonumber\\
 & \simeq & \tilde{\phi}^*i_H^*(\ko_S\boxtimes Q\boxtimes\ko_H)\nonumber\\
 & \simeq & \tilde{\phi}^*(\be_H\otimes L_H).\nonumber
\end{eqnarray}
By functoriality $L_{\tilde{Y}}\simeq\tilde{\phi}^*L_H,$ and hence
$\tilde{\be}\simeq\tilde{\phi}^*\be_H$.
\vskip5mm

The group $SL(N_m)$
acts naturally on the Grassmannian $G=\Grass(N_m,2)$ and induces an action
on $\Hilb^{P_H}(S\times G)$, under which $H$ is invariant. In order to find an
algebraic structure on the quotient by the action of $SL(N_m)$, we have to shrink
$H$. Consider the open $SL(N_m)$-invariant
subsets $H^g$ and $H^t$ of $H$ defined as follows:
$$
H^g:=\{(S_T\overset{i}\hookrightarrow S\times G)
\in H\ |\ (i^*(\ko_S(-mh)\boxtimes Q))/S_T/\C)\in\M_n^g(\C)\},
$$
$$
H^t:=\{(S_T\overset{i}\hookrightarrow S\times G)
\in H\ |\ (i^*(\ko_S(-mh)\boxtimes Q))/S_T/\C)\in\M_n^t(\C)\}.
$$

Note that 
\begin{equation}\label{inters H}
 H^s=H^g\cap H^t.
\end{equation}

It is well known that

(i) the action $SL(N_m)\times H^s\to H^s$ is proper and $H^s//SL(N_m)$ is isomorphic to $M^s_{S,h}(2;\kn,n)$;

(ii) $M_n^t:=H^t//SL(N_m)$ is isomorphic to an open subscheme of the Gieseker-Maruyama moduli space
$M_{S,h}(2;\kn,n)$ containing $M^s_{S,h}(2;\kn,n)$ as a dense open subscheme. 

In fact,
the embeddings into $S\times G$ used in our construction are equivalent to
the embeddings into $G$ used by Gieseker, in the case when the underlying tree-like surface is irreducible, hence the Gieseker (semi)stability condition on the vector bundles from $H^t$ coincides with the Mumford (semi)stability under the action of $SL(N_m)$. This
proves (ii). From Propositions 3.1, 3.2 of \cite{G1}, we conclude that
the points of
$H^t//SL(N_m)$ represent exactly the $S$-equivalence classes of semistable vector bundles on $S$.

Since $M_{S,h}(2;\kn,n)$ is
projective, $M_n^t$ is quasi-projective. In particular, $M_n^t$ is separated.

\begin{sub}\label{properness}{\bf Proposition.}\;
The action
$SL(N_m)\times H^g\to H^g$ is proper.\
\end{sub}

For the proof we use the properness criterion via families over curves and
Theorem \ref{septhm} on separatedness.
\vskip5mm

Now we invoke the following result of Koll\'ar \cite[Theorem 1.5]{Ko1}.

{\bf Theorem.}\;
{\it Fix an excellent base scheme $\Lambda$. Let $G$ be an affine algebraic
group scheme of finite type over $\Lambda$ and $X$ a separated algebraic space of finite type over
$\Lambda$. Let $m:G\times X\to X$ be a proper $G$-action on $X$. Assume that one of the following
conditions is satisfied:

(1)  $G$ is a reductive group scheme over $\Lambda$.

(2)  $\Lambda$ is the spectrum of a field of positive characteristic.

Then a geometric quotient $p_X:X\to X//G$ exists and $X//G$ is a separated algebraic space of
finite type over $\Lambda$.}

Applying now the case (1) of this theorem to our situation with
$\Lambda=\Spec(\C),\ G=SL(N_m),\ X=H^g,$ we obtain from
Proposition \ref{properness}
that $M_n^g:=H^g//SL(N_m)$ is a separated algebraic space of finite type over
$\C$ and
$p:H^g\to H^g//SL(N_m)=M_n^g$ is a geometric quotient.

Note that
$$
M^s_{S,h}(2;\kn,n)=M_n^g\cap M_n^t
$$
is open in both $M_n^g$, $M_n^t$, so we can glue them together along $M_n^s$ 
into an
algebraic space
$$M_n^0=M_n^g\cup M_n^t,$$
which is of finite type, separated but not necessarily complete.

For an arbitrary $\bt_n$-family
$(\be/\bx/Y)$ one has an $SL(N_m)$-equivariant diagram (\ref{univmor}).
Assume in addition that the family belongs to $\M^g(Y)$. Then 
there is the diagram
\begin{equation}\label{univmor2}
\xymatrix{\tilde{Y}\ar[d]_{\rho}\ar[r]^{\phi}& H^g\ar[d]_{p_H}\\
Y & M_n^g .}
\end{equation}
As $\rho$ is a principal bundle map, it follows
that there exists a morphism (of algebraic
spaces) $f:Y\to M_n^g$ which extends (\ref{univmor2}) to a commutative diagram
\begin{equation}\label{univmor3}
\xymatrix{\tilde{Y}\ar[d]_{\rho}\ar[r]^{\phi}& H^g\ar[d]_{p_H}\\
Y\ar[r]^{f}&M_n^g .}
\end{equation}
The existence of this modular morphism $f:Y\to M_n^g$
and the existence of the semi-universal family $(\be_H^g/\bx_H^g/H^g)$ means that 
$M_n^g$
corepresents the functor $\bm_n^g$, i.e. it is the wanted moduli space.
The same argument applies to $M_n^t$ and $M_n^0$.
Thus we have the following

\begin{sub}\label{thmexist}{\bf Theorem.}\;
There exists a separated 
algebraic space $M_n^g$ (resp. $M_n^0$) of finite type corepresenting the functor $\bm_n^g$ (resp. $\bm_n^0$).
\end{sub}

There is a particular case in which we can say more. This is the case when $\M_n^g=\M_n$. This happens when
$M_{S,h}(2;\kn,n)$ contains no strictly semistable locally free sheaves,
for example, when $S=\P_2$. Then we can state:

\begin{sub}\label{thmexist-a}{\bf Theorem.}\; Let $S=\P_2$.
Then $\M_n^g=\M_n$, and there exists a proper\footnote{I. e. complete, separated and of finite type}
algebraic space $M_n$ corepresenting the functor $\bm_n$.
\end{sub}

In the general case, we can only suggest a conjecture.

\begin{sub}\label{thmexist-b}{\bf Conjecture.}\; Let $S$ be any smooth 
projective surface. Then there exists a proper
algebraic space corepresenting $\bm_n$.
\end{sub}

It is not clear whether one should expect $M_n$ to be projective.
However, in some examples one can present an explicit
construction of $M_n$ as a projective variety. One such example is treated
in the next section.

\section{The bubble-tree compactification of $M_{\P_2}(2; 0,2)$}

Let $M(0,2)=M_{\P_2}(2; 0,2)$ be the moduli space of semistable sheaves
with Chern classes $c_1=0, c_2=2$ and rank $2.$ It is well known that
$M(0,2)$ is isomorphic to the $\P_5$ of conics in the dual plane, the
isomorphism being given by $[\kf]\leftrightarrow C(\kf),$ where $C(\kf)$
is the conic of jumping lines of $[\kf]$ in the dual plane, see
\cite{Ba}, \cite{Ma}, \cite{OSS}. For a more explicit description we
use the following notation.

\begin{sub}\label{beilin}{\bf Beilinson resolutions.}\;
\rm In the sequel $S$ will be the
projective plane:
$S=P(V)=\P_2$, where $V$ is a fixed 3-dimensional vector space. We will write 
$m\kf$ for $\C^m\otimes\kf$, where $\kf$ is 
a sheaf. 

It is well known that any sheaf $\kf$ from $M(0,2)$ has
two Beilinson resolutions
$$
0\to 2\, \Omega^2_S(2)
\overset{A}{\lra} 2\, \Omega_S^1(1)\to\kf\to 0
$$
$$0\to 2\, \ko_S(-2)
\overset{B}{\lra} 4\, \ko_S(-1)\to\kf\to 0,
$$
where the matrices $A$ (of vectors in $V$) and $B$ (of vectors in $V^*$)
are related by the exact sequence
$$
0\to \C^2\overset{A}{\lra}\C^2\otimes V \overset{B}{\lra}\C^4\to 0.
$$
Recall that $\Hom(\Omega_S^2(2), \Omega_S^1(1))$ is canonically isomorphic
to $V$ with $v\in V$ acting by contraction. The matrix product $BA$ is zero,
where the elements of the two matrices are multiplied by the rule that for $v\in V$ and $f\in V^*$, the product $fv$ is $f(v).$

The matrices $A$ and $B$ are determined by $\kf$ uniquely up to isomorphisms
of the above resolutions. The first resolution implies that
$\det(A)\in S^2V$ is non-zero. It is the equation of the conic $C(\kf).$
The sheaf $\kf$ is locally free if and only if $C(\kf)$ is smooth, or
if and only if $\kf$ is stable. If $C(\kf)$ decomposes into a pair of lines,
then $A$ is equivalent to a matrix of the form
$\left( \begin{smallmatrix}x& 0\\z & y\end{smallmatrix}\right),$
and $\kf$ is an extension
\[
0\to \ki_{\langle x\rangle }\to \kf\to\ki_{\langle y\rangle }\to 0.
\]
In this case $\kf$ is $S$-equivalent to 
$\ki_{\langle x\rangle }\oplus\ki_{\langle y\rangle }$ .
\end{sub}

We will first present the $\bt_2$-bundles appearing in the compactificaton as
limits of 1-parameter degenerations. We will start from
the following explicit description of the blowup of $\A^1\times S$ at a
point.

\begin{sub}\label{blowZ}{\bf A special blowup.}\;\rm
 Let $(e_0, e_1, e_2)$ be a basis of $V$ and $(x_0, x_1, x_2)$ the 
corresponding homogeneous coordinates on $S$. The blowup
$$X\xra{\sigma}\A^1\times S$$
of the point $p=(0,\langle e_0\rangle)$ is the subvariety of 
$\A^1\times S\times \P_2$
given by the equations
\begin{eqnarray*}
tx_0u_1-x_1u_0 & = & 0\\
tx_0u_2-x_2u_0 & = & 0\\
x_1u_2-x_2u_1 & = & 0,
\end{eqnarray*}
where the $u_\nu$ are the coordinates of the third factor $\P_2.$
We consider the following divisors on $X$:
\begin{enumerate}
\item $\tilde{S}$, the proper transform of $\{0\}\times S$, 
isomorphic to the blowup of $S$ at $p$; 
\item $S_1$, the exceptional divisor of $\sigma$; 
\item $H$, the lift of $\A^1\times h$, where $h$ is a general line in $S$;
\item $L$, the divisor defined by $\ko_X(L)=pr_3^*\ko_{\P_2}(1)$.
\end{enumerate}
We have
$\tilde S\sim H-L$. We also let $x_\nu$ resp. $u_\nu$ denote the sections
of $\ko_X(H)$ resp. $\ko_X(L)$ lifting the above coodinates. Using the
equations of $X$, we see that the canonical
section $s$ of $\ko_X(\tilde S)$ fits into the diagrams
\begin{equation}\label{sects}
\xymatrix{\ko_X\ar[r]^(0.35){s}\ar[dr]_(0.3){(t_0, x_1, x_2)} & \ko_X(H-L)\ar[d]^{(u_0, u_1, u_2)}\\
 & \ko_X(H).}
\end{equation}
\end{sub}


\begin{sub}\label{expl1}{\bf Example.}\;\rm
Using the notation from \ref{blowZ}, we choose
$A=\left( \begin{smallmatrix}e_0 & 0\\ 0 & e_0\end{smallmatrix}\right),$
representing $\kf_0=\ki_0\oplus\ki_0$, where $\ki_0$
is the ideal sheaf of the point $\langle e_0\rangle\in P(V)$. 
Denote $C=\A^1(\C)$
and let $\kf$ be the family of $M(0,2)$-sheaves over $C\times S$ defined
as the cokernel of the matrix
$$A(t)=\left(\begin{matrix}e_0 & -ta\\ -tb& e_0\end{matrix}\right)$$
with $a=\alpha_1e_1+\alpha_2e_2$ and $b=\beta_1e_1+\beta_2e_2.$
Then $\kf|\{t\}\times S$ is locally free for $t\ne 0$. 
The second Beilinson resolution of $\kf$ (see \ref{beilin}) is then 
$$0\to 2\, \ko_C\boxtimes\ko_S(-2)
\overset{B(t)}{\lra} 4\,\ko_C\boxtimes\ko_S(-1)\to\kf\to 0.
$$
with $B(t)$ given by
$$B(t)=\left(\begin{matrix}x_1 & x_2 & \alpha_1tx_0 & \alpha_2tx_0\\
\beta_1tx_0 & \beta_2tx_0 & x_1 & x_2\end{matrix}\right).$$

Let now $X\xra{\sigma}C\times S$ be the blowup of $C\times S$ at
$p=(0,\langle e_0\rangle)$. We have the following commutative
diagram with exact rows and columns:
\[
\xymatrix{ & & & 0\ar[d] &\\
 & 0\ar[d] & & 2\,\ko_{S_1}(-1)\ar[d] & \\
0\ar[r] & 2\,\ko_X(-2H)\ar[r]^{\sigma^*B(t)}\ar[d]^s & 4\,\ko_X(-H)\ar[r]
\ar@{=}[d] & \sigma^*\kf\ar[r]\ar[d] & 0\\
0\ar[r] & 2\,\ko_X(-H-L)\ar[r]^(0.6){B_X}\ar[d] & 4\,\ko_X(-H)\ar[r]& \bbf\ar[r]\ar[d]  & 0\\ & 2\,\ko_{S_1}(-1)\ar[d] & &  0 & \\
 & 0 & & &
}\]

By (\ref{sects}), 
$$B_X=\left(\begin{matrix}u_1 & u_2 & \alpha_1u_0 & \alpha_2u_0\\
\beta_1u_0 & \beta_2u_0 & u_1 & u_2\end{matrix}\right).$$
Here $2\,\ko_{S_1}(-1)$ is the torsion of $\sigma^*\kf$, supported
on the exceptional divisor, and $\bbf$ is its locally free quotient.
Finally, let $\be:=\bbf(-\tilde{S})$. One can verify that the restriction
of $\be$ to $\tilde{S}$ is trivial,
$\be_{\tilde{S}}\simeq 2\,\ko_{\tilde{S}},$ and that the restriction
$\be_{S_1}$ belongs to $M^s_{\P_2}(2;0,2).$
Thus $\be|\tilde{S}\cup S_1$ is a $\bt_2$-bundle with weighted tree of
type (II), see (\ref{figtyp}),
and is a flat degeneration of bundles in $M(0,2)$.
Note that the vectors $a,b$ used in the construction of this degeneration
correspond to two points on the double line $C(\kf_0)=\{e_0^2=0\}$
and thus determine a ``complete conic'' in classical terminology.
Choosing different such pairs will lead to non-isomorphic $\bt_2$-bundles.
\end{sub}

\begin{sub}\label{expl2}{\bf Example.}\;\rm
Let now
$A=\left( \begin{smallmatrix}e_0 & 0\\ 0 & e_2\end{smallmatrix}\right)$
and $\kf_0=\ki_0\oplus\ki_2$, where $\ki_0$ resp.
$\ki_2$
are the ideal sheaves of the points $p_0=\langle e_0\rangle$ resp. 
$p_2=\langle e_2\rangle$ in $S=P(V)$.
Here we consider the blowup
$X\xra{\sigma}C\times S$  of $C\times S$ at the two points $p_0,\, p_2$.
Let $\tilde{S}=S_0$ again denote the proper transform of $\{0\}\times S$,
and let $S_1$ and $S_2$ denote the exceptional divisors of $\sigma$. We embed
$X$ into $C\times S\times\P_2\times\P_2$
with equations
\begin{eqnarray*}
tx_0u_1-x_1u_0 & = & 0\\
tx_0u_2-x_2u_0 & = & 0 \\
x_1u_2-x_2u_1 & = & 0 \\
tx_2v_0-x_0v_2 & = & 0\\
tx_2v_1-x_1v_2 & = & 0\\
x_1v_0-x_0v_1 & = & 0,
\end{eqnarray*}
where $u_\nu$ and $v_\nu$ are the coordinates of the third and fourth
factors $\P_2$ respectively. Define the sheaf $\kf$ over $C\times S$ as the 
cokernel
of the matrix
$$A(t)=\left(\begin{matrix}e_0 & -t^2e_1\\ -t^2e_1& e_2\end{matrix}\right).$$
The corresponding
matrix $B(t)$ is then given by
$$B(t)=\left(\begin{matrix}x_1 & x_2 & 0 & t^2x_0\\
t^2x_2 & 0 & x_0 & x_1\end{matrix}\right).$$
In order to find the limit $\bt_2$-bundle on the $\bt_2$-surface
$\tilde{S}\cup S_1\cup S_2$, we proceed as in Example \ref{expl1}.
Let $x_\nu, u_\nu, v_\nu$ denote the sections
of $\ko_X(H), \ko_X(L_1), \ko_X(L_2)$ obtained by lifting the respective
homogeneous coordinates, and let
$s_1\in\Gamma\ko_X(S_1)$ and $s_2\in\Gamma\ko_X(S_2)$ be the canonical
sections with divisors $S_1\sim H-L_1$ and $S_2\sim H-L_2.$
We obtain a torsion free sheaf $\bbf$ on $X$
as the quotient in the exact triple
$$0\to\ko_{S_1}(-1)\oplus\ko_{S_2}(-1)\to\sigma^*\kf\to\bbf\to 0$$
with resolution
$$0\to\ko_X(-H-L_1)\oplus\ko_X(-H-L_2)\xra{B_X}4\,\ko_X(-H)\to\bbf\to 0,$$
where
 $$B_X=\left(\begin{matrix}u_1 & u_2 & 0 & tu_0\\
tv_2 & 0 & v_0 & v_1\end{matrix}\right).$$
The restrictions of $\bbf$ to $\tilde{S}, S_1, S_2$ are now
$$ \bbf_{\tilde{S}}\simeq\ko_{\tilde{S}}(-\ell_1)\oplus \ko_{\tilde{S}}(-\ell_2),$$
where $\ell_i=\tilde S\cap S_i$ ($i=1,2$) are the two exceptional curves on
$\tilde S$, and
$$\bbf_{S_1}\simeq\ko_{S_1}\oplus \ki_{q_1}(1),\quad
\bbf_{S_2}\simeq\ko_{S_2}\oplus \ki_{q_2}(1),$$
where $q_1$ resp. $q_2,$ are the points $\{u_1=u_2=0\}$ resp.
$\{v_0=v_1=0\}.$

Next, let $\be'$ be the elementary transform given by the exact sequence
$$0\to\be'\to\bbf\to\ko_{S_1}\oplus\ko_{S_2}\to 0.$$
This sheaf turns out to be locally free on $X$. Then the
sheaf $\be:=\be'(-\tilde{S})$ has the restrictions
$$ \be_{\tilde{S}}\simeq 2\,\ko_{\tilde{S}},\quad \be_{S_1},\quad \be_{S_2},$$
where the Chern classes of $\be_{S_1}$, $\be_{S_2}$ are $c_1=0,\, c_2=1$, and
there are non-split extensions of the form
$$0\to\ko_{S_1}\to\be_{S_1}\to\ki_{q_1}\to0,$$
$$0\to\ko_{S_2}\to\be_{S_2}\to\ki_{q_2}\to0.$$
Thus $\be|\tilde{S}\cup S_1\cup S_2$ is a $\bt_2$-bundle with weighted
tree of type (III), see (\ref{figtyp}).
It is a again a flat degeneration of vector bundles from $M(0,2)$.
\end{sub}

The third example with weighted tree of type (IV) can be constructed by a similar
but slightly more complicated procedure. The final result about the moduli space
$M_2$ is the following theorem.

\begin{sub}\label{thm02}{\bf Theorem.}\; Let $\widetilde{P(S^2V)}$ denote the
blowup of $P(S^2V)$ along the Veronese surface.
\begin{itemize}
\item[(1)] The moduli space $M_2$, defined in Theorem \ref{thmexist-a} with $n=2$, is isomorphic
to $\widetilde{P(S^2V)}$.\\
\item[(2)] The isomorphism classes $[E_T,S_T],\; T\in\bt_2,$ are in
1:1 correspondence with the points of $\widetilde{P(S^2V)}$.
\item[(3)] The weighted trees associated to the pairs $[E_T,S_T]$ that occur in $M_2$
are of one of the following four types:
\end{itemize}

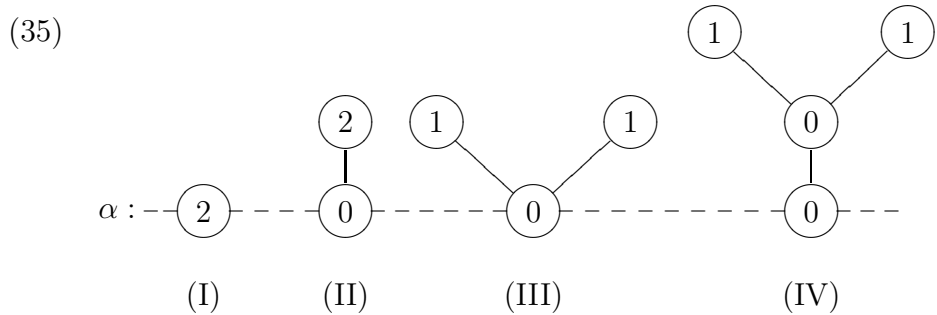
\begin{equation}\label{figtyp}
\xymatrix @-1pc{& & & & & & & *++[o][F-]{1} \ar@{-}[dr]  & & *++[o][F-]{1} \ar@{-}[dl] \\
& & & *++[o][F-]{2} \ar@{-}[d] &  *++[o][F-]{1} \ar@{-}[dr] & & *++[o][F-]{1} \ar@{-}[dl] & &
*++[o][F-]{0} \ar@{-}[d] & \\
\alpha: \ar@{--}[r] & *++[o][F-]{2} \ar@{--}[rr] & & *++[o][F-]{0} \ar@{--}[rr]
& & *++[o][F-]{0} \ar@{--}[rrr] & & & *++[o][F-]{0} \ar@{--}[r] & \\
& {\rm(I)} & & {\rm(II)} & & {\rm(III)}  & & & {\rm(IV)} &
}
\end{equation}

The four types of weighted trees define a stratification of $\widetilde{P(S^2V)}$
in locally closed subsets in the following way. Let $\Sigma_0$ be the exceptional
divisor of $\widetilde{P(S^2V)}$ and $\Sigma_1$ the proper transform
of the cubic hypersurface of decomposable conics in $P(S^2V)$.
Then
\begin{itemize}
\item the points of $\widetilde{P(S^2V)}\smallsetminus\Sigma_0\cup\Sigma_1$
represent the bundles of type (I) on the original surface $S;$
\item the bundles in $\Sigma_0\smallsetminus\Sigma_1$ are of type (II);
\item the bundles in $\Sigma_1\smallsetminus\Sigma_0$ are of type (III);
\item the bundles in $\Sigma_1\cap\Sigma_0$ are of type (IV).
\end{itemize}
\end{sub}

There is a construction of a complete family $(\be/\bx/F)$ of $\bt_2$-bundles
which contains all types of such bundles. This will be sketched next.

\begin{sub}\label{unifam}{\bf Semi-universal family for $M(0,2)$.}\;\rm
From now on $H$ will denote the vector space $\C^2$ and
$G$ will be the Grassmannian $\Grass(2, H\otimes V)=G(2,6)$ of 2-dimenional
subspaces of $H\otimes V$. Let $\ku$ denote the universal subbundle
on $G$. For any subspace $\C^2\overset{y}{\hookrightarrow}H\otimes V$
there is the determinant homomorphism
$\C\intoo{\bigwedge^2y}\wedge^2H\otimes S^2V\simeq S^2V.$
We denote by $G^{ss}$ the open subset defined by  $\wedge^2y\ne 0$.
By the description in \ref{beilin}, this open subset parametrizes all the
sheaves in $M(0,2)$, and there is a semi-universal family ${\boldsymbol{\kf}}$ on 
$G^{ss}\times P$
with resolution
$$0\to\ku\boxtimes\Omega_S^2(2)\to H\otimes\ko_{G^{ss}}\boxtimes\Omega_S^1(1)\to{\boldsymbol{\kf}}\to 0,$$
where $\ku$ is the universal subbundle on $G$, restricted to $G^{ss}$.
The map $y\mapsto\wedge^2y$ defines a morphism $G^{ss}\to P(S^2V)$
which is the modular morphism of the family and at the same time
is a good quotient by
the natural action of $SL(H)$ on $G^{ss}$:
$$G^{ss}//SL(H)\simeq  P(S^2V)\simeq M(0,2).$$
Consider the subvarieties
$\Delta_0,\; \Delta_1^{\prime},\; \Delta_1^{\prime\prime},\; \Delta_1^{\prime\prime\prime}$
of $G^{ss}$ defined as follows.
For each $y\in G^{ss}$, let $l_y\subset P(H\otimes V)$ be the corresponding
line, and let $\mathcal S$ be the image of 
the Segre embedding $P(H)\times P(V)\hookrightarrow P(H\otimes V)$. Then:

\begin{tabular}{lll}
$ \Delta_0$ & $:=$ & $\{y\in G^{ss}|\; l_y\subset \mathcal S\}$,\\
$\Delta_1^{\prime}$ & $:=$ & $\{y\in G^{ss}|\;l_y\cap \mathcal S \;\text{consists of two simple points}\}$,\\
$\Delta_1^{\prime\prime}$ & $:=$ & $\{y\in G^{ss}|\;l_y\cap \mathcal S \;\text{is a double point}\}$,\\
$\Delta_1^{\prime\prime\prime}$ & $:=$ & $\{y\in G^{ss}|\;l_y\cap \mathcal S \;\text{is a simple point}\}$.
\end{tabular}

One finds the following normal forms for the
matrices $A$ defining the inclusions $\C^2\overset{A}{\hookrightarrow}H\otimes V$
that represent the points
$y\in G^{ss}$:
\begin{enumerate}
\item[--]
$y\in \Delta_0$ if and only if $y$ is represented by a matrix
$A=\left( \begin{smallmatrix}x& 0\\0 & x\end{smallmatrix}\right),$
\item[--]
$y\in \Delta_1^{\prime}$ if and only if $y$ is represented by a matrix
$A=\left(\begin{smallmatrix}x& 0\\0 & x^\prime\end{smallmatrix}\right)$
with independent vectors $x$ and $x^\prime,$
\item[--]
$y\in \Delta_1^{\prime\prime}$ if and only if $y$ is represented by a matrix
$A=\left(\begin{smallmatrix}x & 0\\z & x\end{smallmatrix}\right)$
with independent vectors $x$ and $z,$
\item[--]
$y\in \Delta_1^{\prime\prime\prime}$ if and only if $y$ is represented by a matrix
$A=\left(\begin{smallmatrix}x & 0\\z & x^\prime\end{smallmatrix}\right)$,
where ($x$, $x^\prime$, $z$) is a basis ov $V.$
\end{enumerate}
\end{sub}

\begin{sub}\label{kirwan}{\bf The Kirwan blowup.}\;\rm
Let $\tilde{G}\xra{\tilde{\sigma}}G$ be the blowup of $G$ along $\Delta_0$
followed by the blowup along the proper transform
$\tilde{\Delta_1}$ of $\Delta_1:=\Delta_1^{\prime}\cup \Delta_1^{\prime\prime}.$
Then $\tilde{G}$ is also acted on by $SL(H)$ and there is a suitable
linearization of this action such that
$\tilde{G}^s=\tilde{G}^{ss}$.
The geometric quotient $\tilde{G}^s//SL(H)$ is isomorphic to
$\widetilde{P(S^2V)}$, and we have the following commutative diagram:

\begin{equation}\label{kirwqu}
\xymatrix{\tilde{G}^s\ar[d]_{\tilde{\pi}}\ar[r]^{\tilde{\sigma}} & G^{ss}\ar[d]^{\pi}\\
\widetilde{P(S^2V)}\ar[r]^{\sigma} & P(S^2V).}
\end{equation}

We have the following exceptional divisors in $\tilde{G}^s$.
Let $D_0:=\tilde{\sigma}^{-1}(\Delta_0)$, and let $D_1$ be the inverse image
of $\tilde{\Delta_1}$ in $\tilde{G}^s$.
Then
$$
D^s_i:=D_i\cap\tilde{G}^s=\tilde{\pi}^{-1}(\Sigma_i),\ \ i=1,\: 2,$$
where $\Sigma_0$ and
$\Sigma_1$ are the divisors in $\widetilde{P(S^2V)}$ introduced in Theorem~\ref{thm02}.
\end{sub}

\begin{sub}\label{bigfam}{\bf Univesal family of Serre constructions.}\;\rm
We parametrized the first Beilinson resolutions
$$
0\to 2\ko_S(-1)\to 2\Omega^1_S(1)\to\kf\to 0,
$$
of the sheaves from $M(0,2)$ by the two-dimensional
subspaces $U$ of $H\otimes V=H^0(S,H\otimes\Omega^1_S(2))$.
Now we want to parametrize all the Serre constructions
of these sheaves. A Serre construction for $\kf$ is determined
by a section $s$ of $\kf(1)$ with finitely many zeros. 
The global sections of $\kf(1)$, up to proportionality, are in one-to-one
correspondence with three dimensional subspaces $W$ of $H\otimes V$
such that $U\subset W$. Such a subspace provides an extension of
the injection in the Beilinson resolution
$2\ko_S(-1)=U\otimes\ko_S(-1)\hookrightarrow 2\Omega^1_S(1)=H\otimes\Omega^1_S(1)$ to
an injection $W\otimes\ko_S(-1)\hookrightarrow H\otimes\Omega_S(1)$, determined
up to the action of $GL(W)$, and this extension defines
two exact sequences
$$
0\to W\otimes\ko_S(-1)\to H\otimes\Omega^1_S(1)\to \ki_{Z,S}(1)\to 0, 
$$
$$
0\to(W/U)\otimes \ko_S(-1)\xrightarrow{s}\kf\to \ki_{Z,S}(1)\to 0.
$$
Here $Z$ is a zero-dimensional subscheme of $S$ of length 3.
The second sequence represents the Serre construction for $\kf$. 
For a given $\kf$, we can always find $W$ such that $Z$ is either the union 
of three distinct points, or an isotropic fat point 
Spec$\:\ko_S/\mathfrak m_x^2$ ($x\in S$) of length 3.
Let now relativize this construction over the whole of $\tilde{G}^s$.
Let
$$F_{2,3}:=\{(U,W)|\; U\subset W\subset H\otimes V\},$$
be the flag variety of 2- and 3-dimensional subspaces of $H\otimes V$, and
$$F=F_{2,3}\times_G\tilde{G}^s .$$ There are
natural projections
$\tilde{G}^s\overset{\gamma}\leftarrow F\overset{\delta}\to F_{2,3}$.
We have the semi-universal family of Beilinson resolutions of the polystable 
sheaves from $M(0,2)$
over $\tilde G^s$. Lifting it to $F$, we obtain the exact sequence 
$$
0\to\bu\boxtimes\Omega^2_S(2)\to
H\otimes\ko_F\boxtimes\Omega^1_S(1)\to{\boldsymbol{\kf}}\to0 .
$$
Shrinking $F$ to an appropriate open subset, mapped surjectively
onto $\tilde G^s$, we obtain
a semi-universal family of Serre constructions over $F$,
\begin{equation}\label{univ_Serre}
0\to\bw/\bu\boxtimes\Omega^2_S(2)\to{\boldsymbol{\kf}}\to\ki_{\kz,F\times S}\otimes
\ko_F\boxtimes\ko_S(1)\to 0,
\end{equation}
where $\bu$ and $\bw$ are the lifts of the tautological subbundles
from $F_{2,3}$, and where $\kz$ is a flat family of zero-dimensional 
subschemes of $S$
of length 3 over $F$. We can shrink further $F$ in such a way that $\gamma$
remains surjective, so that the only singularities
of $\kz$ are the quasi-transversal intersections of three smooth branches 
over the points of $D_0$.
\end{sub}
\vskip5mm

\begin{sub}\label{finblow}{\bf Semi-universal family for $M_2$.}\;\rm
Let $\D_0:=\gamma^{-1}(D_0^s)$ and $\D_1:=\gamma^{-1}(D_1^s)$ be the
lifted divisors in $F$, and let
$$B_0:=\kz\cap (\D_0\times S)$$
be the codimension 3 intersection. This is the singular locus 
of $\kz$, where three branches intersect. 
Besides $B_0$, the singular locus $\Sing {\boldsymbol{\kf}}$ of ${\boldsymbol{\kf}}$ contains
the points $(f,x_i)\in \kz$ such that $f\in \D_1$, $x_i\in S$ ($i=1,2$), 
$x_1\neq x_2$ and  
${\boldsymbol{\kf}}|_{f\times S}\simeq \ki_{x_1,S}\oplus \ki_{x_2,S}$. 
At these points
$\kz$ is smooth, but the local extension class of the Serre sequence 
degenerates. We are to resolve both types of singularities.

Let
$$\bx^\prime\xra{\sigma_0}F\times S$$
be the blowup of $B_0.$ Let further $B_1$ be the closure of
$$\sigma_0^{-1}(\Sing {\boldsymbol{\kf}}\cup ((F-\D_0)\times S)).$$
Let 
$$\bx\xra{\sigma_1}\bx^\prime$$
be the blowup of $B_1$. Consider the composed morphism
$$p: \bx\xrightarrow{\sigma_1\circ\sigma_0}F\times 
S\overset{\raisebox{0.2ex}{$\scriptstyle{pr}$}_{\scriptscriptstyle 2}}{\lra}F,$$
and let $\E_0$, $\E_1$ denote the exceptional divisors of the last two blowups.
Then $\bx\xra{p}F$ is a family of $\bt_2$-surfaces which
includes the above examples. Moreover, the proper transform $\tilde\kz$ of 
$\kz$
in $\bx$ is smooth. In order to replace it by a vector bundle, we consider
the Serre construction on $\bx$ lifting \eqref{univ_Serre}.
A computation shows that the local extension
class of \eqref{univ_Serre}, when lifted to a local section of the invertible sheaf
$\epsilon\in \mathcal Ext^1\big(\ki_{\tilde Z,\bx}\otimes pr_S^*\ko_S(1), p^*(\bw/\bu\boxtimes\Omega^2_S(2))\big)$,
has a simple pole along $\D_0\cap\kz$ and is regular and nonvanishing everywhere else.

To transform the simple pole into a double one, we make the base change $\tilde{\bx}
\to \bx$, which is a double covering branched at $\D_0+\D_1$. It is defined locally
over $\bx$. We add $\ \tilde{}\ $ to mark the lifts to $\tilde{\bx}$ of all the objects defined
on $\bx$. Then $\tilde\epsilon$ acquires a double pole along $\tilde{\D}_0
\cap\Tilde{\Tilde{\kz}}$ and has no
other singularities. Hence it defines a regular nowhere vanishing
section of the invertible sheaf
$\mathcal Ext^1\big(\ki_{\Tilde{\Tilde{\kz}},\tilde{\bx}}\otimes pr_S^*\ko_S(1), \tilde p^*(\bw/\bu\boxtimes\Omega^2_S(2))(2\tilde{\D}_0)\big)$. This implies that in the extension
$$
0\to \tilde p^*(\bw/\bu\boxtimes\Omega^2_S(2))(\tilde{\D}_0) \to
\be \to 
\ki_{\Tilde{\Tilde{\kz}},\tilde{\bx}}\otimes pr_S^*\ko_S(1)(-\tilde{\D}_0)
\to 0
$$
defined by $\tilde\epsilon$, the middle term $\be$ is a vector bundle.
This is the wanted family of $\bt_2$-bundles. It is defined over $\tilde{\bx}$,
which can be thought of as a DM stack with stabilizers of order $\leq 2$
whose associated coarse moduli space is $\bx$. The statements of Theorem \ref{thm02}
follow by considering the classifying map of our functor $\bm_2$ on this family towards
the moduli space $M_2$.
\end{sub}



\begin{thebibliography}{MMMMM}

\bibitem[Ba]{Ba} {\bf Barth W.} {\it Moduli of vector bundles on the
projective plane,} Invent. math. 42 (1977), 63--91.


\bibitem[B1]{B1} {\bf Buchdahl N. P.} {\em Blowups and gauge fields},  Pacific J. Math.  {\bf 196}, 69--111  (2000).

\bibitem[B2]{B2} {\bf Buchdahl N. P.} {\em Sequences of stable bundles over compact complex surfaces},  J. Geom. Anal. {\bf 9}, 391--428  (1999).




\bibitem[D1]{D1} {\bf Donaldson S. K.} {\em Anti self-dual Yang-Mills connections over complex algebraic surfaces and stable vector bundles}, Proc. London Math. Soc. (3) 50 (1985), no. 1, 1--26.

\bibitem[D2]{D} {\bf Donaldson S. K.} {\em Compactification and completion of Yang-Mills moduli spaces},  Differential geometry (Pe\~n\'\i scola, 1988),  145--160, Lecture Notes in Math., 1410, Springer, Berlin, 1989.





\bibitem[Fe]{Fe} {\bf Feehan P.M.N.} {\it Geometry of the ends of the moduli space of
anti-self-dual connections,} J.Diff.G. {\bf42}, No.3 (1995), 465--553.



\bibitem[FM]{FM} {\bf Fulton W., MacPherson R.} {\it A Compactification of
configuration spaces,} Annals of Math. 139 (1994), 183--225.

\bibitem[G1]{G1} {\bf Gieseker D.} {\it On the moduli of vector bundles on an algebraic surface},  Ann. of Math. (2)  {\bf 106}  (1977), 45--60.

\bibitem[G2]{G} {\bf Gieseker D.} {\it A construction of stable bundles on
an algebraic surface,} J. Differential Geom. 27 (1988), 137--154.







\bibitem[HL]{HL} {\bf Huybrechts D., Lehn M.} {\it The Geometry of Moduli
Spaces of Sheaves,} Aspects of Mathematics E 31, Braunschweig: Vieweg 1997.



\bibitem[Ki]{Ki} {\bf Kirwan F.} {\it Partial desingularisations of quotients
of nonsingular varieties and their Betti numbers,}
Ann of Math. 122 (1985), 41--85.


\bibitem[Ko]{Ko} {\bf Koll\'ar J}. {\it Projectivity of complete moduli,}
J. Differential Geometry, 32 (1990), 235--268.

\bibitem[Ko1]{Ko1} {\bf Koll\'ar J.} {\it Quotient spaces modulo algebraic
groups,} Annals of Math. 145 (1997), 33--79.



\bibitem[JL]{JL} {\bf Li, Jun}: {\it  Algebraic geometric interpretation of
Donaldson's polynomial invariants,} J. Diff. Geom. 37, 417--466, 1993

\bibitem[LT]{LT} {\bf  L\"ubke M., Teleman, A.}: {\it  The Kobayashi-Hitchin correspondence}, World Scientific, River Edge, 1995.

\bibitem[Ma]{Ma} {\bf Maruyama M.} {\it Singularities of the curves of jumping
lines of a vector bundle of rank 2 on $\mathbb{P}_2$.} Algebraic Geometry,
Proc.of Japan-France Conf.,
Tokyo and Kyoto, 1982, Lect. Notes in Math., 1016, Springer, 1983, 370--411.

\bibitem[MT1]{MT1} {\bf Maruyama M., Trautmann G.}
{\it On compactifications of the moduli space of instantons},
Intern. Journal of Math. 1, (1990), 431--477.

\bibitem[MT2]{MT2} {\bf Maruyama M., Trautmann G.}:
{\it Limits of instantons}, Intern. Journal of Math. 3, (1992), 213--276.




\bibitem[NS]{NS} {\bf Nagaraj D. S., Seshadri C. S.}
{\em Degenerations of the moduli spaces of vector bundles on curves},
I,  Proc. Indian Acad. Sci. Math. Sci. {\bf 107}, 101--137  (1997);
II, Proc. Indian Acad. Sci. Math. Sci. {\bf 109}, 165--201  (1999).


\bibitem[OSS]{OSS}{\bf Okonek C., Schneider M., Spindler H.}
{\it Vector Bundles on Complex Projective Spaces}, Birkh\"auser, 1980.








\bibitem[Sim]{Sim} {\bf Simpson C.T.} {\it Moduli of representations of
the fundamental group of a smooth projective variety I,}
Publ. Math. I.H.E.S. 79 (1994), 47--129.

\bibitem[Ta]{Ta1} {\bf Taubes C.H.} {\it A framework for Morse theory for the Yang-Mills
functional,} Invent. math. {\bf94} (1988), 327--402.


\bibitem[Tr]{T} {\bf Trautmann G.} {\it Moduli spaces in algebraic geometry,}
Manuscript, 2000.

\bibitem[V]{V} {\bf Viehweg E.} {\it Quasi-projective moduli for polarized
manifolds,} Springer, 1995.

\bibitem[U]{U} {\bf Uhlenbeck K. K.} {\em Removable singularities in Yang--Mills fields},  Comm. Math. Phys. {\bf 83}, 11--29  (1982).


\end{thebibliography}
\end{document}